\documentclass[notitlepage,leqno,10pt]{article}
\usepackage{latexsym}
\usepackage{amsmath}
\usepackage{amsfonts}
\usepackage{amssymb}

\renewcommand{\d}{\delta }
\newcommand{\D }{\Delta }

\newcommand{\e }{\varepsilon }

\renewcommand{\l }{\lambda }

\newcommand{\n }{\nabla }
\newcommand{\var }{\varphi }

\newcommand{\rh }{\rho }

\newcommand{\s }{\sigma }
\newcommand{\Sig }{\Sigma}

\newcommand{\ov}{\overline}
\newcommand{\intbar}{\mathop{\int\makebox(-13.5,0){\rule[4pt]{.7em}{0.3pt}}%
\kern-6pt}\nolimits}

\newcommand{\be}{\begin{equation}}
\newcommand{\ee}{\end{equation}}
\newcommand{\bes}{\begin{equation*}}
\newcommand{\ees}{\end{equation*}}
\newcommand{\ba}{\begin{eqnarray}}
\newcommand{\ea}{\end{eqnarray}}
\newcommand{\bas}{\begin{eqnarray*}}
\newcommand{\eas}{\end{eqnarray*}}
\newenvironment{pf}{\noindent{\sc Proof}.\enspace}{\rule{2mm}{2mm}\medskip}
\newenvironment{pfn}{\noindent{\sc Proof}}{\rule{2mm}{2mm}\medskip}

\newtheorem{remark}{Remark}[section]
\newcommand{\R}{\mathbb{R}}

\newcommand{\Z}{\mathbb{Z}}

\newcommand{\N}{\mathbb{N}}
\def\dint{\displaystyle{\int}}
\author{ Cheikh Birahim NDIAYE}

\date{}

\title{\bf Constant $T$-curvature conformal metrics on \;$4$-manifolds with boundary}
\begin{document}

\newtheorem{lem}{Lemma}[section]
\newtheorem{pro}[lem]{Proposition}
\newtheorem{thm}[lem]{Theorem}
\newtheorem{rem}[lem]{Remark}
\newtheorem{cor}[lem]{Corollary}
\newtheorem{df}[lem]{Definition}

\maketitle

\begin{center}

{\small  SISSA, via Beirut 2-4, 34014 Trieste, Italy.}

\end{center}

\

\

\begin{center} 
{\bf Abstract}
\end{center}
In this paper we prove that, given a compact four dimensional smooth Riemannian manifold \;$(M,g)$\;with smooth boundary there exists a metric conformal to\;$g$\;with constant \;$T$-curvature, zero\;$Q$-curvature and zero mean curvature under generic and conformally invariant assumptions. The problem amounts to solving a fourth order nonlinear elliptic boundary value problem (BVP) with boundary conditions given by a third-order pseudodifferential operator, and homogeneous Neumann one. It has a variational structure, but since the corresponding Euler-Lagrange functional is in general unbounded from below, we look for saddle points. In order to do this, we use topological arguments and  min-max methods combined with a compactness result for the corresponding BVP.

\begin{center}

\bigskip\bigskip
\noindent{\bf Key Words:} Geometric BVPs, Blow-up analysis, Variational methods, Min-max schemes\\ \;$Q$-curvature, \;$T$-curvature, Conformal geometry, Topological methods.
\bigskip

\centerline{\bf AMS subject classification: 35B33, 35J35, 53A30, 53C21}

\end{center}
\section{Introduction}
In recent years, there has been an intensive study of conformally covariant differential (or even pseudodifferential) operators on compact smooth Riemannian manifolds, their associated curvature invariants in order to understand the relationships between analytic and geometric properties of such objects.
\vspace{2pt}

\indent
A model example is the Laplace-Beltrami operator on compact closed surfaces \;($\Sigma ,g$), which governs the transformation laws of the Gauss curvature. In fact under  the conformal change of metric \;$ g_u=e^{2u}g$,\; we have 
\ba\label{eq:g}
\Delta_{g_u}=e^{-2u}\Delta_{g};\;\;\;\;\;\;\;\;\;-\Delta_{g}u+K_{g}=K_{g_u}e^{2u},
\ea
where \;$\Delta_{g}$\;and\;$K_{g}$\;(resp.\;$\Delta_{g_u}$\;and\;$K_{g_u}$)\; are the Laplace-Beltrami operator and the Gauss curvature of ($\Sigma,g$) (resp. of ($\Sigma_,g_u$)).\\
Moreover we have the Gauss-Bonnet formula which relates  \;$\int_{\Sigma}K_{g}dV_{g}$\;and the topology of \;$\Sigma$\;:
$$
\int_{\Sigma}K_{g}dV_{g}=2\pi\chi(\Sigma);
$$ 
where\;$\chi(\Sigma)$\; is the Euler-Poincar\'e  characteristic of \;$\Sigma$. From this we deduce that \;$\int_{\Sigma}K_{g}dV_{g}$\;is a topological invariant (hence also a conformal one). Of particular interest is the classical {\em Uniformization Theorem} which says that every compact closed Riemannian surface carries a conformal metric with constant Gauss curvature.

There exists also a conformally covariant differential operator on four dimensional compact closed Riemannian manifolds called the Paneitz operator, and to which is associated a natural concept of curvature. This operator, discovered by Paneitz in 1983 (see\;$\cite{p1}$)\;and the corresponding Q-curvature  introduced by Branson (see\;$\cite{bo}$)\; are defined in terms of Ricci tensor \;$Ric_{g}$\; and scalar curvature \;$R_{g}$\; of the Riemannian manifold \;$(M,g)$\; as follows
\ba\label{eq:P}
P_{g}\varphi=\D_{g}^{2}\varphi+div_{g}(\frac{2}{3}R_{g}g-2Ric_{g})d\varphi;
\ea
\begin{equation}\label{eq:qcur}
Q_{g}=-\frac{1}{12}(\D_{g}R_{g}-R_{g}^{2}+3|Ric_{g}|^{2}),
\end{equation}
where\;$\varphi$\; is any smooth function on \;$M$.\\

\indent
As the Laplace-Beltrami operator governs the transformation laws of the Gauss curvature, we also have  that the Paneitz operator does the same for the Q-curvature. Indeed under a conformal change of metric \;$g_u=e^{2u}g$\;we have  
$$\label{eq:tranlaw}
P_{g_u}=e^{-4u}P_{g};\;\;\;\;\;\;\;\;\;P_{g}u+2Q_{g}=2Q_{g_u}e^{4u}.
$$
\indent
Apart from this analogy, we also have an extension of the Gauss-Bonnet formula which is the Gauss-Bonnet-Chern formula
$$
\int_{M}(Q_{g}+\frac{|W_{g}|^{2}}{8})dV_{g}=4\pi^{2}\chi(M),
$$
where \;$W_{g}$\;denotes the Weyl tensor of \;($M,g$), see\;$\cite{dm}$.\;Hence,  from the pointwise conformal invariance of\;$|W_{g}|^{2}dV_{g}$, it follows that the integral of \;$Q_{g}$\;over \;$M$\;is also a conformal invariant one.\\

\indent
As for the {\em Uniformization Theorem} for compact closed Riemannian surfaces, one can also ask if every closed compact four dimensional Riemannian manifolds carries a metric conformally related to the background one with constant \;$Q$-curvature.\\
A first positive answer to this question was given by Chang-Yang\cite{cy} under the assumptions that \;$P_g$\;non-negative and \;$\int_{M}Q_gdV_g<8\pi^2$. Later Djadli-Malchiodi\cite{dm} extend Chang-Yang result to a large class of compact closed four dimensional Riemannian manifold assuming that \;$P_g$\;has no kernel and \;$\int_{M}Q_gdV_g$\;is not an integer multiple of \;$8\pi^2$.\\

\indent
On the other hand, there are high-order analogues to the Laplace-Beltrami operator and to the Paneitz operator for high-dimensional compact closed Riemannian manifolds and also to the associated curvatures (called again \;$Q$-curvatures), see\cite{fg},\cite{fg1} and \cite{gjms}.\\
As for the question of the existence of constant \;$Q$-curvature conformal metrics on a given compact closed four dimensional Riemannian manifold,  regarding high-dimensional \;$Q$-curvature, one can still ask the same question for a compact closed Riemannian manifolds of arbitrary dimensions.\\
A first affirmative answer has been given by Brendle in the {\em even} dimensional case under the assumption that the high-dimensional analogue of the Paneitz operator is non-negative and the total integral of the \;$Q$-curvature is less than \;$(n-1)!\omega_n$ ( where \;$\omega_n$\;is the area of the unit sphere \;$S^n$\;of \;\;$R^{n+1}$)\;using a geometric flow, see\cite{bren}. The result of Djadli-Malchiodi\cite{dm} (and the one in\cite{bren}) has been extended to all dimensions in \cite{nd}. \\

\noindent
As for the case of compact closed Riemannian manifolds, many works have also been done in the study of conformally covariant differential operators on compact smooth Riemannian manifolds with smooth boundary, their associated curvature invariants, the corresponding boundary operators and curvatures in order also to understand the relationship between analytic and geometric properties of such objects.\\

\indent
A model example is the Laplace-Beltrami operator on compact smooth surfaces  with smooth boundary\;($\Sigma ,g$),  and the Neumann operator on the boundary. Under a conformal change of metric the couple constituted by the Laplace-Beltrami operator and the Neumann operator govern the transformation laws of the Gauss curvature and the geodesic curvature. In fact, under the conformal change of metric \;$g_u=e^{2u}g$,\;we have 
\begin{equation*}
\left\{
     \begin{split}
\D_{g_u}=e^{-2u}\D_g;\\
\frac{\partial }{\partial n_{g_u}}=e^{-u}\frac{\partial }{\partial n_{g}};
     \end{split}
     \right.
\qquad \mbox{and} \qquad
\left\{
     \begin{split}
      -\D_{g}u+K_{g}=K_{g_u}e^{2u}\;\;\text{in}\;\;\Sigma;\\
        \frac{\partial u}{\partial n_{g}}+k_g=k_{g_u}e^{u}\;\;\;\text{on}\;\;\partial \Sigma.
     \end{split}
   \right.
\end{equation*}
where \;$\Delta_{g}$\;(resp. \;$\D_{g_u}$)\;is the Laplace-Beltrami operator of (\;$\Sigma ,g$)\;(resp. ($\Sigma_,g_u$))\; and\;$K_{g}$\;(resp. \;$K_{ g_u}$)\;is the Gauss curvature of ($\Sigma,g$) (resp. of ($\Sigma_,g_u$)), \;$\frac{\partial }{\partial n_{g}}$\;(resp \;$\frac{\partial }{\partial n_{g_u}}$)\;is the Neumann operator of (\;$\Sigma ,g$)\;(resp. of ($\Sigma_,g_u$)) and \;$k_{g}$\; (resp.\;$k_{g_u}$) is the geodesic curvature of \;($\partial \Sigma,g$) (resp of ($\partial \Sigma_, g_u$)) .\\
Moreover we have the Gauss-Bonnet formula which relates  \;$\int_{\Sigma}K_{g}dV_{g}+\int_{\partial \Sigma}k_gdS_g$\;and the topology of \;$\Sigma$
\begin{equation}\label{eq:gb}
\int_{\Sigma}K_{g}dV_{g}+\int_{\partial \Sigma}k_gdS_g=2\pi\chi(\Sigma),
\end{equation}
where\;$\chi(\Sigma)$\; is the Euler-Poincar\'e  characteristic of \;$\Sigma$,\;$dV_g$\;is the element area of \;$\Sig$ and \;$dS_g$\;is the line element of \;$\partial \Sig$. Thus\;$ \int_{\Sigma}K_{g}dV_{g}+\int_{\partial \Sigma}k_gdS_g$\;is a topological invariant, hence a conformal one.\\
 In this context, of particular interest is also an analogue of the  classical {\em Uniformization Theorem}, namely given a compact Riemannian surface\;$(\Sig,g)$\; with boundary, does there exists metrics conformally related to\;$g$\;with constant Gauss curvature and constant geodesic curvature. This problem has been solved through the following theorem (for a proof see \cite{br})
\begin{thm}\label{eq:unif}
Every compact smooth Riemannian surface with smooth boundary\;$(\Sigma,g)$\;carries a metric conformally related to \;$g$\;with constant
Gauss curvature and constant geodesic curvature.
\end{thm}

\vspace{2pt}

\indent
As for compact closed four dimensional Riemannian manifolds, on four-manifolds with boundary we also have  the Paneitz operator \;$P^4_g$\; and the \;$Q$-curvature. They are defined with the same formulas (see \;$\eqref{eq:P}$\;and\;$\eqref{eq:qcur}$) and enjoy the same invariance properties as in the case without boundary, see\;$\eqref{eq:tranlaw}$.\\

\indent
Likewise, Chang and Qing\cite{cq1} have  discovered a boundary operator \;$P^3_g$\;defined on the boundary of compact four dimensional smooth Riemannian manifolds and a natural third-order curvature \;$T_g$\;associated to \;$P^3_g$ as follows
\begin{equation*}
P^3_g\varphi=\frac{1}{2}\frac{\partial {\D_g\varphi}}{\partial n_g}+\D_{\hat g}\frac{\partial \varphi}{\partial n_g}-2H_g\D_{\hat g}\varphi+(L_g)_{ab}(\nabla_{\hat g})_{a}(\nabla_{\hat g})_{b}+\nabla_{\hat g}H_g.\nabla_{\hat g}\varphi+(F-\frac{R_g}{3})\frac{\partial \varphi}{\partial n_g}.
\end{equation*}
\begin{equation*}
T_g=-\frac{1}{12}\frac{\partial R_g}{\partial n_g}+\frac{1}{2}R_gH_g-<G_g,L_g>+3H_g^3-\frac{1}{3}Tr(L^3)+\D_{\hat g}H_g,
\end{equation*} 
where\;$\varphi$\; is any smooth function on \;$M$,\;\;\;$\hat g$\;is the metric induced by \;$g$\;on\;\;$\partial M$, \;$L_g=(L_g)_{ab}=-\frac{1}{2}\frac{\partial g_{ab}}{\partial n_g}$\\is the second fundamental form of \;$\partial M$,\;\;$H_g=\frac{1}{3}tr(L_g)=\frac{1}{3}g^{ab}L_{ab}$\;(\;$g^{a,b}$\;are the entries of the inverse \;$g^{-1}$\;of the metric\; $g$)\;is the mean curvature of \;$\partial M$,\;$R^k_{bcd}$\;is the Riemann curvature tensor\;\; $F=R^{a}_{nan}$,\;\;$R_{abcd}=g_{ak}R^{k}_{bcd}$\;(\;$g_{a,k}$\;are the entries of the metric \;$g$)\;and\;\;$<G_g,L_g>=R_{anbn}(L_g)_{ab}$.\\

\indent
On the other hand, as the Laplace-Beltrami operator and the Neumann operator govern the transformation laws of the Gauss curvature and the geodesic curvature on compact surfaces with boundary under conformal change of metrics, we have that the couple\;$(P^4_g,P^3_g)$\;does the same for \;$(Q_g,T_g)$\;on compact four dimensional smooth Riemannian manifolds with smooth boundary. In fact, after a conformal change of metric \;$ g_u=e^{2u}g$\;we have that
\begin{equation*}
\left\{
     \begin{split}
P^4_{g_u}=e^{-4u}P^4_g;\\
P^3_{g_u}=e^{-3u}P^3_{g}; 
     \end{split}
   \right.
\qquad \mbox{and}\qquad
\left\{
\begin{split}
P^4_g+2Q_g=2Q_{ g_u}e^{4u}\;\;\text{in }\;\;M\\
P^3_g+T_g=T_{ g_u}e^{3u}\;\;\text{on}\;\;\partial M.
\end{split}
\right.
\end{equation*}
\indent
Apart from this analogy we have also an extension of the Gauss-Bonnet formula \;$\eqref{eq:gb}$\; which is known as the Gauss-Bonnet-Chern formula
\begin{equation}\label{eq:gbc}
\int_{M}(Q_{g}+\frac{|W_{g}|^{2}}{8})dV_{g}+\int_{\partial M}(T+Z)dS_g=4\pi^{2}\chi(M)
\end{equation}
where \;$W_g$\;denote the Weyl tensor of \;$(M,g)$\; and \;$ZdS_g$\;(for the definition of \;$Z$\;see \cite{cq1}) are pointwise conformally invariant. Moreover, it turns out that \;$Z$\;vanishes when the boundary is totally geodesic (by totally geodesic we mean that the boundary \;$\partial M$\;is umbilic and minimal).\\
Setting
\begin{equation*}
\kappa_{P^4_g}=\int_{M}Q_gdV_g,\;\;\;\;\;\kappa_{P^3_g}=\int_{\partial M}T_gdS_{g};
\end{equation*}
we have that thanks to \;$\eqref{eq:gbc}$, and to the fact that \;$W_gdV_g$\;and \;$ZdS_g$\;are pointwise conformally invariant,\;$\kappa_{P^4_g}+\kappa_{P^3_g}$\;is conformally invariant, and will be denoted by
\begin{equation}\label{eq:invc}
\kappa_{(P^4,P^3)}=\kappa_{P^4_g}+\kappa_{P^3_g}.
\end{equation}

\indent
The Riemann mapping Theorem is one of the most celebrated theorems in mathematics. It says that an open, simply connected, proper subset of the plane is conformally diffeomorphic to the disk. So one can ask if such a theorem remains true in dimension \;$4$. Unfortunately in dimension 4 few regions are conformally diffeomorphic to the ball. \\
However, in the spirit of the Uniformization Theorem  (Theorem\;$\ref{eq:unif}$), one can still ask, if on a given compact four dimensional smooth Riemannian manifold with smooth boundary, there exists a metric conformal to the background one with zero\;$Q$-curvature, constant \;$T$-curvature  and  zero mean curvature.\\
In the context of the Yamabe problem, related questions were raised by Escobar \cite{es}.\\
\vspace{2pt}

\indent
In this paper, we are interested to give an analogue of the Riemann mapping Theorem (in the spirit of Theorem\;$\ref{eq:unif}$) to compact four dimensional smooth  Riemannian manifold with smooth boundary under generic and conformally invariant assumptions.
Writting \;$ g_u=e^{2u}g$, the problem is equivalent to solving the following BVP:
\begin{equation*}
\left\{
\begin{split}
P^4_gu+2Q_g&=0\;\;&\text{in}\;\;M;\\
P^3_gu+T_g&=\bar Te^{3u}\;\;&\text{on}\;\;\partial M;\\
\frac{\partial u}{\partial n_g}-H_gu&=0\;\;&\text{on}\;\;\partial M.
\end{split}
\right.
\end{equation*}
where \;$\bar Q$\;is a fixed real number and\;$\frac{\partial}{\partial n_g}$\;is the inward normal derivative with respect to \;$g$. 
\vspace{4pt}

\indent
Due to a result by Escobar, \cite{es}, and to the fact that we are interested to solve the problem under conformally invariant assumptions, it is not restrictive to assume\;$H_g=0$, since this can be always obtained through a conformal transformation of the background metric. Thus we are lead to solve the following BVP with Neumann homogeneous boundary condition:
\begin{equation}\label{eq:bvps}
\left\{
\begin{split}
P^4_gu+2Q_g&=0\;\;&\text{in}\;\;M;\\
P^3_gu+T_g&=\bar Te^{3u}\;\;&\text{on}\;\;\partial M;\\
\frac{\partial u}{\partial n_g}&=0\;\;&\text{on}\;\;\partial M.
\end{split}
\right.
\end{equation}
\vspace{2pt}

\noindent
Defining \;$H_{\frac{\partial}{\partial n}}$\;as
\begin{equation*}
H_{\frac{\partial}{\partial n}}=\Big\{u\in H^2(M):\;\;\;\frac{\partial u}{\partial n_g}=0\Big\};
\end{equation*}
and \;$P^{4,3}_g$\;as follows, for every \;$u,v\in H_{\frac{\partial }{\partial n}}$
\begin{equation*}
\begin{split}
\left<P^{4,3}_gu,v\right>_{L^2(M)}=\int_{M}\left(\D_g u\D_gv+\frac{2}{3}R_g\nabla_g u\nabla_g v\right)dV_g-2\int_{M}Ric_g(\nabla_g u,\nabla_g v)dV_g\\-2\int_{\partial M}L_g( \nabla_{\hat g} u, \nabla_{\hat g} v)dS_g,
\end{split}
\end{equation*}
we have that by the regularity result in Proposition\;$\ref{eq:reg}$\;\;below, critical points of the functional
\begin{equation*}
\begin{split}
II(u)=\left<P^{4,3}u,u\right>_{L^2(M)}+4\int_{M}Q_gudV_g+4\int_{\partial M}T_gudS_g-
\frac{4}{3}\kappa_{(P^4,P^3)}\log\int_{\partial M}e^{3u}dS_g;\;\;\;u\in H_{\frac{\partial }{\partial n}},
\end{split}
\end{equation*}
which are weak solutions of \;$\eqref{eq:bvps}$\; are also smooth and hence strong solutions.\\
\vspace{2pt}

\indent
A similar problem has been adressed in \cite{nd1}, where constant \;$Q$-curvature metrics with zero \;$T$-curvature and zero mean curvature are found under generic and conformally invariant assumptions.\\\\
In \cite{nd2}, using heat flow methods, it is proven that if the operator\;$P^{4,3}_g$\;is non-negative, $Ker P^{4,3}_g\simeq \R$, and \;$\kappa_{(P^4,P^3)}< 4\pi^2$ the problem \;$\eqref{eq:bvps}$ is solvable.\\
Here we are interested to extend the above result under generic and conformally invariant assumptions.\\
Our main theorem is:
\begin{thm}\label{eq:theorem3}
Suppose  \;$Ker P^{4,3}_g\simeq\R$. Then assuming\;$\kappa_{(P^4,P^3)}\neq k4\pi^2$\;for\;$k=1,2,\cdots$, we have that \;$(M,g)$\;admits a conformal metric with constant\;$T$-curvature, zero\;$Q$-curvature and zero mean curvature.
\end{thm}
\begin{rem}
a) Our assumptions are conformally invariant and generic, so the result applies to a large class of compact\;$4$-dimensional manifolds with boundary.\\
b) From the Gauss-Bonnet-Chern formula, see\;$\eqref{eq:gbc}$\; we have that Theorem\;$\ref{eq:theorem3}$\;does NOT cover the case of locally conformally flat manifolds with totally geodesic boundary and positive integer Euler-Poincar\'e characteristic.\\
\end{rem}
Our assumptions include the two following situations:\\
\begin{equation}\label{eq:sit1}
\kappa_{(P^{4},P^3)}<4\pi^2 \;\;\text{and (or)}\;\;P^{4,3}_{g}\;\text{possesses}\;\;\bar k\;\;\text{ negative eigenvalues (counted with multiplicity)}
\end{equation}
\begin{equation}\label{eq:range}
\begin{split}
\kappa_{(P^{4},P^3)}\in \left(4k\pi^2 \;\;, \;\;4(k+1)\pi^2\right),\;\;\;\text{for some }\;\;k\in \N^*  \;\;\text{and (or)}\;\;P^{4,3}_{g}\;\text{possesses}\;\;\bar k\;\;\text{negative eigenvalues}\\\text{ (counted with multiplicity)}
\end{split}
\end{equation}
\begin{rem}
Case \;$\eqref{eq:sit1}$\;includes the condition ($\bar k=0$)\;under which in \cite{nd2} it is proven existence of solutions to \;$\eqref{eq:bvps}$, hence will not be considered here. However due to a trace Moser-Trudinger type inequality (see Proposition\;$\ref{eq:mos-trub}$\;below) it can be achieved using  Direct Method of Calculus of Variations. \\
In order to simplify the exposition, we will give the proof of Theorem\;$\ref{eq:theorem3}$\;in the case where we are in situation \;$\eqref{eq:range}$\; and \;$\bar k=0$\;(namely \;$P^{4,3}_g$\;is non-negative). At the end of Section 4 a discussion to settle the general case \;$\eqref{eq:range}$\; and also case\;$\eqref{eq:sit1}$\; is made.
\end{rem}
\vspace{2pt}

\noindent
To prove Theorem\;$\ref{eq:theorem3}$\;we look for critical points of\;$II$. Unless \;$\kappa_{(P^4,P^3)}<4\pi^2$\; and \;$\bar k=0$, this Euler-Lagrange functional is unbounded from above and  below (see Section 4), so it is necessary to find extremals which are possibly saddle points. To do this we will use a min-max method: by classical arguments in critical point theory, the scheme yields a {\em Palais-Smale sequence}, namely a sequence \;$(u_l)_l\in H_{\frac{\partial }{\partial n}}$\;satisfying the following properties
\begin{equation*}
II(u_l)\rightarrow c\in\R;\;\;\;\;II ^{'}(u_l)\rightarrow 0\;\;\text{as}\;l\rightarrow+\infty.
\end{equation*}
Then, as is usually done in min-max theory, to recover existence one should prove that the so-called {\em Palais-Smale condition} holds, namely that every Palais-Smale sequence has a converging subsequense or a similar compactness criterion. Since we do not know if the Palais-Smale condition holds, we will employ Struwe's monotonicity method, see \cite{str}, also used in \cite{dm} and \cite{nd}. The latter yields existence of solutions for arbitrary small perturbations of the given equation, so to consider the original problem one is lead to study compactness of solutions to perturbations of \;$\eqref{eq:bvps}$. Precisely we consider \\
\begin{equation}\label{eq:sep}
\left\{
\begin{split}
P^4_gu_l+2Q_l&=0\;\;&\text{in}\;\;M;\\
P^3_gu_l+T_l&=\bar T_le^{3u_l}\;\;&\text{on}\;\;\partial M;\\
\frac{\partial u_l}{\partial n_g}&=0\;\;&\text{on}\;\;\partial M.
\end{split}
\right.
\end{equation}
where
\begin{equation}\label{eq:lbq}
\bar T_l \longrightarrow \bar T_{0}>0\;\;\;\;in\;\;C^{2}(\partial M)\;\;\;T_{l}\longrightarrow  T_{0}\;\;\;\;in\;\;C^{2}(\partial M)\;\;\;\;\;Q_{l}\longrightarrow  Q_{0}\;\;\;\;in\;\;C^{2}(M);\;
\end{equation}
\vspace{4pt}

\noindent
\begin{rem}\label{eq:blowbo}
From the Green representation formula given in Lemma\;$\ref{eq:greenrep}$ below, we have that if \;$u_l$\;is a sequence of solutions to \;$\eqref{eq:sep}$, then \;$u_l$\;satisfies 
$$
u_l(x)=-2\int_{M}G(x,y)Q_l(y)dV_g-2\int_{\partial M}G(x,y)T_l(y)dS_g(y)+2\int_{\partial M}G(x,y)\bar T_l(y)e^{3u_l(y)}dS_g(y).
$$
Therefore, under the assumption\;$\eqref{eq:lbq}$, if\;$\sup_{\partial M}u_l\leq C$, then we have \;$u_l$\;is bounded in \;$C^{4+\alpha}$\;for every \;$\alpha\in(0,1)$.
\end{rem}

\indent
In this context, due to Remark\;$\ref{eq:blowbo}$\; we say that a sequence\;$(u_l)$\;of solutions to \;$\eqref{eq:sep}$\;{\em blows up} if the  following holds:
\begin{equation}\label{eq:blowup}
\text{there exist}\;\; x_l\in \partial M\;\;\text{such that} \;\;u_l(x_l)\rightarrow+\infty\;\;\text{as}\;\;l\rightarrow+\infty,
\end{equation}
and we prove the following compactness result.
\begin{thm}\label{eq:comp}
Suppose\;$Ker P^{4,3}_g\simeq\R$\;and that \;$(u_{l})$\;is a sequence of solutions to\;$\eqref{eq:sep}$\;with\;\;$\bar T_{l}$,\;$T_{l}$\; and \;$Q_l$\;satisfying\;$\eqref{eq:lbq}$. Assuming that
\;$(u_l)_l$\;blows up (in the sense of \;$\eqref{eq:blowup}$) and 
\ba\label{eq:limee}
\int_{M}Q_0dV_g+\int_{\partial M}T_0dS_g+o_l(1)=\int_{\partial M}\bar 
T_le^{3u_l}dS_g;
\ea 
then there exists \;$N\in\N\setminus\{0\}$\;such that 
\begin{equation*}
\int_{M}Q_{0}dV_{g}+\int_{\partial M}T_{0}dS_g=4N\pi^2.
\end{equation*}
\end{thm}
\vspace{6pt}

\noindent
From this we derive a corollary which will be used to ensure compactness of some solutions to a sequence of approximate BVP's produced by the topological argument combined with Struwe's monotonicity method. Its proof is a trivial application of Theorem\;$\ref{eq:comp}$\;and\; Proposition\;$\ref{eq:reg}$\;below.
\begin{cor}\label{eq:compco}
 Suppose\;$Ker P^{4,3}_{g}\simeq\R$.\\
 a) Let\;$(u_{l})$\;be a sequence of solutions to\;$\eqref{eq:sep}$\;with\;\;$\bar T_{l}$,\;$T_{l}$\; and \;$Q_l$\;satisfying\;$\eqref{eq:lbq}$. Assume also that
 \begin{equation*}
\int_{M}Q_0dV_g+\int_{\partial M}T_0dS_g+o_l(1)=\int_{\partial M}\bar T_le^{3u_l}dV_g;
\end{equation*} 
 and 
\begin {equation*}
k_{0}=\int_{M}Q_{0}dV_{g}+\int_{\partial M}T_0dS_g\neq 4k\pi^2\;\;k=1,2,3,\dots.
\end{equation*}
then \;$(u_l)_l$\;is bounded in \;$C^{4+\alpha}(M)$\;for any \;$\alpha\in (0,1)$.\\\\
b) Let $(u_{l})$\;be a sequence of solutions to\;$\eqref{eq:bvps}$\;for a fixed value of the constant\;$\bar T$. Assume also that\;$\kappa_{(P^4,P^3)}\neq 4k\pi^2$, then \;$(u_l)_l$\;is bounded in \;$C^{m}(M)$\;for every positive integer\;$m$.\\\\
c) Let $(u_{\rho_k})$\;\;$\rho_k\rightarrow 1$\;be a family of solutions to\;$\eqref{eq:bvps}$\;with \;$T_g$\; replaced by \;$\rho_k T_g$\;, \;$Q_g$\;by \;$\rho_k Q_g$\;and \;$\bar T$\;by \;$\rho_k\bar T$\;for a fixed value of the constant\;$\bar T$. Assume also that\;$\kappa_{(P^4,P^3)}\neq 4k\pi^2$, then \;$(u_{\rho_k})_{k}$\;is bounded in \;$C^{m}(M)$\;for every positive integer\;$m$.\\\\
d) If $\kappa_{(P^4,P^3)}\neq 4k\pi^2\;\;k=1,2,3,\dots$, then the set of metrics conformal to\;$g$\;with constant \;$T$-curvature the constant being the same for all of them, and  with zero\;$Q$-curvature and zero mean curvature  is compact in \;$C^m(M)$\;for \;positive integer \;$m$.\\\\
f) If $\kappa_{(P^4,P^3)}\neq 4k\pi^2\;\;k=1,2,3,\dots$, then the set of metrics conformal to\;$g$\;with constant \;$T$-curvature ,zero \;$Q$-curvature, zero mean curvature and of unit boundary volume is compact in \;$C^m(M)$\;for every positive integer \;$m$.
\end{cor}
\vspace{6pt}

\noindent

We are going to describe the main ideas to prove the above results. Since the proof of Theorem\;$\ref{eq:theorem3}$\;relies on the compactness result of Theorem\;$\ref{eq:comp}$\;(see corollary\;$\ref{eq:compco}$), it is convenient to discuss first the latter. We use the same arguments as in \cite{nd} and \cite{nd1} and noticing that, due to the Green representation formula (see Lemma\;$\ref{eq:greenrep}$), we have only to take care of the behaviour of the restriction of \;$u_l$\;on \;$\partial M$, see Remark\;$\ref{eq:blowbo}$.
\vspace{4pt}

\indent
Now having this compactness result we can describe the proof of Theorem\;$\ref{eq:theorem3}$\;assuming \;$\eqref{eq:range}$\;and that\;$P^{4,3}_g$\;is non-negative.  First of all from  \;$\kappa_{(P^4,P^3)}\in (k4\pi^2,(k+1)4\pi^2)$\; and considerations coming from an improvement of Moser-Trudinger inequality, it follows that if \;$II(u)$\;attains large negative values then \;$e^{3u}$ has to concentrate near at most \;$k$\;points of\;$\partial M$. This means that, if we normalize \;$u$\;so that \;\;$\int_{\partial M}e^{3u}ds_g=1$, then naively \;$e^{3u}\simeq\sum_{i=1}^kt_i\d_{x_i},\;\;x_i\in \partial M,\;\; t_i\geq 0,\;\sum_{i=1}^kt_i=1$. Such a family of convex combination of Dirac deltas are called formal barycenters of \;$\partial M$\;of order \;$k$, see Section 2 , and will be denoted by \;$\partial M_k$. With a further analysis (see Proposition\;$\eqref{eq:pro2}$ ), it is possible to show that the sublevel\;$\{II<-L\}$\;for large \;$L$\; has the same homology as \;$\partial M_k$. Using  the non contractibility of \;$\partial M_k$, we define a min-max scheme for a perturbed functional \;$II_{\rho}$, \;$\rho$\;close to \;$1$, finding a P-S sequence to some levels \;$c_{\rho}$. Applying the monotonicity procedure of Struwe, we can show existence of critical points of \;$II_{\rho}$\;for a.e \;$\rho$, and we reduce ourselves to the assumptions of Theorem\;$\ref{eq:compco}$.
\vspace{10pt}

\noindent
{\bf Acknowledgements:}
I would like to thank Professor Andrea Malchiodi for several stimulating discussions. \\
The author have been supported by M.U.R.S.T within the PRIN 2006 Variational methods and nonlinear differential equations.

\section{Notation and Preliminaries}
In this brief section we collect some useful notations, state a lemma giving the existence of the Green function of the operator \;$(P^{4}_{g},P^3_g)$\; with its asymptotics near the singularity and a trace analogue of the well-known Moser-Trudinger inequality for the operator \;$P^{4,3}_{g}$\;when it is non-negative.\\

In the following \;$B_{p}(r)$\; stands for the metric ball of radius \;$r$\;and center\;$p$ ,\;$B^+_{p}(r)=B_{p}(r)\cap M$\;if \;$p\in \partial M$. Sometimes we use\; $B^+_{p}(r)$\;to denote\;$B_{p}(r)\cap M$\;even if \;$p\notin  \partial M$.\\ In the sequel,\;$B^x(r)$\; will stand for the Euclidean ball of center \;$x$\; and radius \;$r$, $B^x_+(r)=B^x(r)\cap \R^4_+$\;if\;$x\in\partial \R^4_+$. We use also\;$B^x_+(r)$\;to denote \;$B^x(r)\cap \R^4_+$\;even if \;$x\notin \partial \R^4_+$. We denote by \;$d_{g}(x,y)$\; the metric distance between two points \;$x$\;and \;$y$\; of \;$M$\, and \;$d_{\hat g}(x,y)$\;the intrinsic distance of two points \;$x$\;and \;$y$\;of \;$\partial M$. Given a point \;$x\in \partial M$, and \;$r>0$,\;$B^{\partial M}_x(r)$\; stands for the metric ball in \;$\partial M$\;with respect to the (intrinsic) distance \;$d_{\hat g}(\cdot,\cdot)$\;of center \;$x$\;and radius \;$r$.\; \;$H^{2}(M)$\;stands for the usual Sobolev space of functions on \;$M$\;which are of class \;$H^2$\;in each coordinate system.  Large positive constants are always denoted by \;$C$,\;and the value of\;$C$\;is allowed to vary from formula to formula and also within the same line.\;$M^2$\;stands for the Cartesian product \;$M\times M$, while \;$Diag(M)$\; is the diagonal of \;$M^2$. Given a function \;$u\in L^1(\partial M)$,\;$\bar u_{\partial M}$\; denotes its average on \;$\partial M$,\;that is \;$\bar u_{\partial M}=\left( Vol_{\hat g}(\partial M)\right)^{-1}\int_{\partial M} u(x)dS_{g}(x)$\; where \;$Vol_{\hat g}(\partial M)=\int_{\partial M}dS_{g}$.\\
$\N$\;denotes the set of non-negative integers.\\
$\N^*$\;stands for the set of positive integers.\\
$A_{l}=o_{l}(1)$\; means that \;$A_{l}\longrightarrow 0$\;\;as the integer \;\;$l\longrightarrow +\infty$.\\
$A_{\epsilon}=o_{\epsilon}(1)$\; means that \;$A_{\epsilon}\longrightarrow 0$\;\;as the real number\;\;$\epsilon\longrightarrow 0$.\\
$A_{\delta}=o_{\delta}(1)$\; means that \;$A_{\delta}\longrightarrow 0$\;\;as the real number \;\;$\delta\longrightarrow 0$.\\
$A_{l}=O(B_{l})$\; means that\;$A_{l}\leq C B_{l}$\;\;for some fixed constant \;$C$.. \\
$dV_g$\;denotes the Riemannian measure associated to the metric\;$g$.\\
$dS_g$\;stands for the Riemannian measure associated to the metric\;$\hat g$\;induced by \;$g$\;on \;$\partial M$.\\
$d\s_{\hat g}$\;stands for the surface measure on boundary of balls of \;$\partial M$.\\
$|\cdot|_{\hat g}$\;stands for the norm associated to \;$g$.\\
$f=f(a,b,c,...)$\;means that \;$f$\;is a quantity which depends only on \;$a,b,c,...$.\\
Next we let \;$\partial M_k$ denotes the family of formal sums
\begin{equation}\label{eq:base}
\partial M_k
=\lbrace \sum_{i=1}^kt_i\d_{x_i},\;\;t_i\geq 0,\;\;\sum_{i=1}^kt_i=1; x_i\in \partial M\rbrace,
\end{equation}
It is known in the literature as the formal set of barycenters relative to \;$\partial M$\;of order \;$k$. We recall that\; $\partial M_k$\;is a stratified set namely a union of sets of different dimension with maximum one equal to \;$4k-1$.\\

Next we recall the following result (see Lemma 3.7 in \cite{dm}), which is necessary in order to carry
out the topological argument below.

\begin{lem}\label{l:nonco} (well-known)
For any $k \geq 1$ one has $H_{4k-1}(\partial M_k;\Z_2) \neq 0$. As a
consequence $\partial M_k$ is non-contractible.
\end{lem}

\noindent If $\var \in C^1(\partial M)$ and if $\s \in \partial M_k$, we denote
the action of $\s$ on $\var$ as
$$
  \langle \s, \var \rangle = \sum_{i=1}^k t_i \var(x_i), \qquad \quad
\s = \sum_{i=1}^k t_i \d_{x_i}.
$$
Moreover, if $f$ is a non-negative $L^1$ function on $\partial M$ with
$\int_{\partial M} f ds_g = 1$, we can define a distance of $f$ from
$\partial M_k$ in the following way
\begin{equation}\label{eq:distf}
    d(f, \partial M_k) = \inf_{\s \in \partial M_k} \sup \left\{ \left|
  \int_{\partial M} f \var dS_g - \langle \s, \var \rangle \right|
  \; | \; \|\var\|_{C^1(\partial M)} = 1 \right\}.
\end{equation}
We also let
\begin{equation}\label{eq:dm}
  \mathcal{D}_{\e,k} = \left\{ f \in L^1(\partial M) \; : \; f \geq 0,
  \|f\|_{L^1(\partial M)} = 1, d(f, \partial M_k) < \e \right\}.
\end{equation}

Now we state a Lemma which asserts the existence of the Green function of \;$(P^4_g,P^3_g)$\;with homogeneous Neumann condition. Its proof can be found in \cite{nd}.
\begin{lem}\label{eq:greenrep}
Assume that \;$Ker P^{4,3}_g\simeq \R$, then the Green function \;$G(x,y)$\;of \;$(P^4_g,P^3_g)$\; exists in the following sense :\\
a) For all functions \;$u\in C^{2}(M),\;\;\frac{\partial u}{\partial n_g}=0$, we have 
\begin{equation*}
u(x)-\bar u=\int_{M}G(x,y)P^{4}_{g}u(y)dV_{g}(y)+2\int_{\partial M}G(x,y')P^3_gu(t)dS_g(y')\;\;\;\;\;\;x\in M
\end{equation*}
b)
\begin{equation*}
G(x,y)=H(x,y)+K(x,y)
\end{equation*}
is smooth on \;$M^2\setminus Diag(M^2)$,\;$K$\;extends to a \;$C^{2+\alpha}$\;function on \;$M^2$ \;and
\begin{equation*}
H(x,y)=
  \left\{
    \begin{array}{ll}
      \frac{1}{8\pi^2}f(r)\log\frac{1}{r}\;\;\;\;\;\text{if}\;\;\;\;B_{\delta}(x)\cap\partial M=\emptyset;&\\\\
      \frac{1}{8\pi^2}f(r)(\log\frac{1}{r}+\log\frac{1}{\bar r})\;\;\;\;\text{otherwise}.&
  \end{array}
  \right .
\end{equation*}
where\;$f(\cdot)=1$\;in\;$[-\frac{\delta}{2},\frac{\delta}{2}]$\;\;and \;$f(\cdot)\in C^{\infty}_0(-\delta,\delta)$, \;$\delta\leq\frac{1}{2}\min\{\delta_1,\delta_2\}$,\;$\delta_1$\;is the injectivity radius of \;$M$\;in \;$\tilde M$, and \;$\delta_2=\frac{\delta_0}{2}$, \;$r=d_g(x,y)$\;and\;$\bar r=d_g(x,\bar y)$.
\end{lem}
\vspace{6pt}

\indent
Next we give a regularity result corresponding to boundary value problems of the type of  BVP\;$\eqref{eq:bvps}$ \;and high order {\em a priori} estimates for sequences of solutions to BVP like\;$\eqref{eq:sep}$\;when they are bounded from above. Its proof is a trivial adaptation of the arguments of Proposition 2.3 in \cite{nd1} 
\begin{lem}\label{eq:reg}
Let \;$u\in H_{\frac{\partial }{\partial n}}$\;be a weak solution to
\begin{equation*}
\left\{
\begin{split}
P^4_gu&=h\;\;&\text{in}\;\;M;\\
P^3_gu+f&=\bar fe^{3u}\;\;&\text{on}\;\;\partial M.
\end{split}
\right.
\end{equation*}
with\;$f\in C^{\infty}(\partial M)$,\;$h\in C^{\infty }(M)$\;and \;$\bar f$\;a real constant. Then we have that \;$u\in C^{\infty}(M)$.\\
Let \;$u_l\in H_{\frac{\partial }{\partial n}}$\;be a sequence of weak solutions to
\begin{equation*}
\left\{
\begin{split}
P^4_gu_l&=h_l\;\;&\text{in}\;\;M;\\
P^3_gu_l+f_l&=\bar f_le^{3u_l}\;\;&\text{on}\;\;\partial M.
\end{split}
\right.
\end{equation*}
with\;$f_l \rightarrow f_0\;\text{in}\; C^{k}(\partial M)$,\;$\bar f_l\rightarrow \bar f_0\;\text{in}\; C^{k}(\partial M)$\;and \; $h_l\rightarrow h_0\;\text{in}\;C^{k}( M)$\;for some fixed \;$k\in \N^*$. Assuming \;$\sup_{\partial M}u_l\leq C$\; we have that
\begin{equation*}
||u_l||_{C^{k+3+\alpha}(M)}\leq C
\end{equation*}
for any \;$\alpha\in(0,1)$.
\end{lem}

\vspace{2pt}

\indent
Now we give a Proposition which is a trace  Moser-Trudinger type inequality when the operator \;$P^{4,3}_g$\;is non-negative with trivial kernel. Its proof can be found in \cite{nd2}, but for the reader convenience we will repeat it here.
\begin{pro}\label{eq:mos-trub}
Assume \;$P^{4,3}_{g}$\;is a non-negative operator with \;$Ker P^{4,3}_{g}\simeq\R$. Then  we have that for all \;$\alpha <12\pi^2$\;there exists a constant \;$C=C(M,g,\alpha)$\;such that 
\begin{equation}\label{eq:mts}
\dint_{\partial M}e^{\frac{\alpha(u-\bar u_{\partial M})^2}{\left<P^{4,3}_{g}u,u\right>_{L^2(M,g)}}}dS_{g}\leq C,
\end{equation}
for all\;$u\in H_{\frac{\partial }{\partial n}}$,\;and hence
\begin{equation}\label{eq:mts}
\log\int_{\partial M}e^{3(u-\bar u)}dS_{g}\leq C+\frac{9}{4\alpha}\left<P^{4,3}_{g}u,u\right>_{L^2(M,g)}\;\;\forall u\in H_{\frac{\partial }{\partial n}}.
\end{equation}
\end{pro}
\begin{pf}
First of all, without loss of generality we can assume \;$\bar u_{\partial M}=0$. Following the same argument as in Lemma 2.2 in \cite{cq2}. we get \;$\forall \beta<16\pi^2$\;there exists \;$C=C(\beta,M)$
\begin{equation*}
\dint_{M}e^{\frac{\beta v^2}{\int_{M}|\D_g v|^2dV_g}}dV_{g}\leq C,\;\;\forall v\in H_{\frac{\partial }{\partial n}}\;\;\text{with}\;\;\bar v_{\partial M}=0.
\end{equation*}
From this, using the same reasoning as in Proposition 2.7 in \cite{nd1}, we derive
\begin{equation}\label{eq:mts}
\dint_{M}e^{\frac{\beta v^2}{\left<P^{4,3}_{g}v,v\right>_{L^2(M)}}}dV_{g}\leq C,\;\;\forall v\in H_{\frac{\partial }{\partial n}}\;\;\text{with}\;\; \bar v_{\partial M}=0.
\end{equation}
Now let \;$X$\;be a vector field extending the the outward normal at the boundary \;$\partial M$. Using the divergence theorem we obtain
$$
\int_{\partial M}e^{\alpha u^2}dS_{g}=\int_{M}div_g\left(Xe^{\alpha u^2}\right)dV_{g}.
$$
Using the formula for the divergence of the product of a vector fied and a function we get
\begin{equation}\label{eq:imp00}
\int_{\partial M}e^{\alpha u^2}dS_{g}=\int_{M}\left(div_gX+2u\alpha \n_{g}u\n_{g}X \right)e^{\alpha u^2}dV_{g}.
\end{equation}
Now we suppose \;$<P^{4,3}_{g}u,u>_{L^2(M)}\leq 1$, then since the vector field \;$X$\;is smooth we have
\begin{equation}\label{eq:imp0}
\left|\int_{M}div_gXe^{\alpha u^2}dV_{g}\right|\leq C;
\end{equation}
thansk to \;$\eqref{eq:mts}$.
Next let us show that 
\begin{equation*}
\left|\int_{M}2\alpha u\n_{g}u\n_{g}X e^{\alpha u^2}dV_{g}.
\right|\leq C
\end{equation*}
Let \;$\epsilon>0$\;small and let us set \;$$p_1=\frac{4}{3-\epsilon},\;\;p_2=4,\;\;p_3=\frac{4}{\epsilon}.$$
It is easy to check that 
$$\frac{1}{p_1}+\frac{1}{p_2}+\frac{1}{p_3}=1.$$
Using Young's inequality we obtain
$$
\left|\int_{M}2\alpha u\n_{g}u\n_{g}X e^{\alpha u^2}dV_{g}
\right|\leq C||u||_{L^{\frac{4}{\epsilon}}}||\n_{g}u||_{L^4}\left(\int_{M}e^{\alpha\frac{4}{3-\epsilon}u^2}dV_{g}\right)^{\frac{3-\epsilon}{4}}.
$$
On the other hand, Lemma 2.8 in \cite{nd1} and Sobolev embedding theorem imply 
$$
||u||_{L^{\frac{4}{\epsilon}}}\leq C;
$$
and
$$
||\n_{g}u||_{L^4}\leq C.
$$
Furthermore from the fact that \;$\alpha<12\pi^2$, by taking \;$\epsilon$\;sufficiently small and using \;$\eqref{eq:mts}$, we obtain
$$
\left(\int_{M}e^{\alpha\frac{4}{3-\epsilon}u^2}dV_{g}\right)^{\frac{3-\epsilon}{4}}.
$$
Thus we arrive to
\begin{equation}\label{eq:imp1}
\left|\int_{M}2\alpha u\n_{g}u\n_{g}X e^{\alpha u^2}dV_{g}
\right|\leq C.
\end{equation}
Hence \;$\eqref{eq:imp00}$,\;$\eqref{eq:imp0}$\;and \;$\eqref{eq:imp1}$\;imply
$$
\int_{\partial M}e^{\alpha u^2}dS_{g}\leq C,
$$
as desired. So the first point of the Lemma is proved.\\
Now using the algebraic inequality 
$$3ab\leq 3\gamma^2a^2+\frac{3b^2}{4\gamma^2},$$
we have that the second point follows directly from the first one. Hence the Lemma is proved.
\end{pf}

\section{Proof of Theorem\;$\ref{eq:comp}$}
This section is concerned about the proof of Theorem\;$\ref{eq:comp}$. We use  the same strategy as in \cite{nd} and \cite{nd1}. Hence in many steps we will be sketchy and referring to the corresponding arguments in \cite{nd}. However, in contrast to the situation in \cite{nd1}, due remark\;$\ref{eq:blowbo}$, we have only to take care of the behaviour of the restriction of the sequence \;$u_l$\;to the boundary\;$M$. 
\vspace{6pt}

\noindent
\begin{pfn}{ of Theorem\;$\ref{eq:comp}$}
\vspace{4pt}

\noindent
First of all, we recall the following particular case of the result of X. Xu ( Theorem 1.2 in \;$\cite{xu}$).
\begin{thm}\label{eq:xu}($\cite{xu}$)
There exists a dimensional constant \;$\s_3>0$\;such that, if\;$u\in C^{1}(\R^3)$\;is solution of the integral equation 
\begin{equation}\nonumber
u(x)=\int_{\R^3}\s_3\log\left(\frac{|y|}{|x-y|}\right)e^{3u(y)}dy +c_0,
\end{equation}
where \;$c_0$\;is a real number, then\;$e^{u}\in  L^{3}(\R^3)$ implies,\;there exists  \;$\l>0$\;and \;$x_0\in \R^3$\;such that 
\begin{equation}\nonumber
u(x)=\log\left(\frac{2\l}{\l^2+|x-x_0|^2}\right).
\end{equation}
\end{thm}
Now, if  \;$\s_3$\;in\;Theorem\;$\ref{eq:xu}$\; we set 
 \;$k_3=2\pi^2\s_3$\;and\;$\gamma_3=2(k_{3})^{3}$

We divide the proof in \;$5$-steps as in \cite{nd}.
\vspace{4pt}

\noindent
{\bf Step 1}\\
\vspace{1pt}

\noindent
There exists \;$N\in \N^*$,\;\;$N$\;converging points \;\; $(x_{i,l})\subset \partial M$\;\; $i=1,...,N $,\;\;$N$\;with limit points \;$x_i\in \partial M$, sequences $(\mu_{i,l}) \;\;i=1;...;N$\;;\;of positive real numbers  converging to \;$0$\;\;such that the following hold: \\\\
$a$)
\\
$$
\hspace{-130pt}\frac{d_{g}(x_{i,l},x_{j,l})}{\mu_{i,l}}\longrightarrow +\infty \;\;\;\; i\neq j \;\;i,j =1,..,N\;\;\;\; and \;\;\;\;\bar T_{l}(x_{i,l})\mu_{i.l}^{3}e^{3u_{l}(x_{i,l})}=1;
$$
$b$)For every \;$i$\\
$$
v_{i,l}(x)=u_{l}(exp_{x_{i,l}}(\mu_{i,l}x))-u_{l}(x_{i,l})-\frac{1}{3}\log(k_3)\longrightarrow V_{0}(x)\;\;\;\;\;\; in\; \;\;C^{1}_{loc}(\R^4_+),\;{V_0}_{|\partial \R^4_+}(x):=\log(\frac{4\gamma_{3}}{4\gamma_{3}^{2}+|x|^{2}});
$$
and
\\ 
$$\hspace{-110pt}
\lim_{R\rightarrow +\infty}\lim_{l\rightarrow +\infty}\int_{B^+_{x_{i,l}}(R\mu_{i,l})\cap \partial M}\bar T_{l}(y)e^{3u_{l}(y)}ds_{g}(y)=4\pi^2;
$$
$c$)\\
\\
$$\hspace{-100pt}There\; exists\;\;C>0\;\; such \;that\; 
\inf_{i=1,...,N}d_{g}(x_{i,l},x)^{3}e^{3u_{l}(x)}\leq C \;\;\;\;\forall x\in \partial M,\;\;\forall l\in \N.
$$

\vspace{10pt}

\noindent
{\em Proof of Step 1}
\vspace{6pt}

\noindent
First of all let\;$x_{l}\in \partial M$\;be such that \;$u_{l}(x_{l})=\max_{x\in \partial M}u_{l}(x)$, then using the fact that \;$u_l$\;blows up we infer\; $u_{l}(x_{l})\longrightarrow +\infty.$\\
Now since \;$\partial M$\;is compact, without loss of generality we can assume that \;$x_l\rightarrow \bar x\in \partial M$.\\
Next let \;$\mu_{l}>0$\; be such that \;$\bar T_{l}( x_l)\mu_{l}^{3}e^{3u_{l}(x_{l})}=1$.\; Since \;$\bar T_{l}\longrightarrow \bar T_{0}\;\; C^{1}(\partial M)$,\; $\bar T_{0}>0$\;and\; $u_{l}(x_{l})\longrightarrow +\infty$, \;we have that \;$\mu_{l}\longrightarrow 0.$\\  
Let \;$ B^{0}_+(\d\mu_{l}^{-1})$\;be the half Euclidean ball of center \;$0$\; and radius\; $\d\mu^{-1}_{l}$,\;with\;$\d>0$\; small fixed . For \;$x\in B^{0}_+(\d\mu_{l}^{-1})$,\; we set
\ba\label{eq:4n}
v_{l}(x)=u_{l}(exp_{x_{l}}(\mu_{l}x))-u_{l}(x_{l})-\frac{1}{3}\log(k_3);
\ea
\ba\label{eq:5n}
\tilde Q_{l}(x)=Q_{l}(exp_{x_{l}}(\mu_{l}x));
\ea
\ba\label{eq:6n}
\tilde {\bar {Q_{l}}}(x)=\bar Q_{l}(exp_{x_{l}}(\mu_{l}x));
\ea
\ba\label{eq:7n}
 g_{l}(x)=\left( exp^{*}_{x_{l}}g\right) (\mu_{l}x).
\ea
Now from the Green representation formula we have, 
\ba\label{eq:tdivn}
u_{l}(x)-\bar u_{l}=\int_{M}G(x,y)P^{4}_{g}u_{l}(y)dV_{g}(y)+2\int_{\partial M}G(x,y')P^3_gu_l(y')dS_g(y');\;\;\; \forall x \in M,
\ea
where \;$G$\;is the Green function of \;$(P^{4}_{g},P^3_g)$\;(see Lemma\;$\ref{eq:greenrep}$).\\
Now using equation\;$\eqref{eq:sep}$\; and differentiating\;$\eqref{eq:tdivn}$\;with respect to \;$x$\;we obtain that for \;$k=1,2$ 
$$
|\n^k u_{l}|_{g}(x)\leq\int_{\partial M}|\n^k G(x,y)|_{g}\bar T_{l}(y)e^{3u_{l}(y)}dV_{g} +O(1),\\
$$
since $ T_{l}\longrightarrow  T_{0}$\; in \;$C^{1}(\partial M)$\;and\;$Q_l\rightarrow Q_0$\;in \;$C^1(M)$.\\
Now let \;$y_{l} \in B^+_{x_{l}}(R\mu_{l}),$\;$R>0$\; fixed, by using the same argument as in \cite{nd}( formula 43 page 11) we obtain
\ba\label{eq:0bad1n}
\int_{\partial M}|\n^k G(y_{l},y)|_{g}e^{3u_{l}(y)}dV_{g}(y)=O(\mu_{l}^{-k})
\ea
Hence we get
\begin{equation}\label{eq:derbvn}
|\n^k v_{l}|_{g}(x)\leq C.
\end{equation}
Furthermore from the definition of\;$v_{l}$\;(see \;$\eqref{eq:4n}$), we get
\begin{equation}\label{eq:ineqvn}
v_{l}(x)\leq v_{l}(0)=-\frac{1}{3}\log(k_3)\;\;\forall x\in \R^{4}_+
\end{equation}
Thus we infer that\;$(v_{l})_l$\; is uniformly bounded in \;$C^{2}(K)$\; for all compact subsets \;$K$\;of \;$\R^{4}_+$. Hence by Arzel\`a-Ascoli theorem we derive that
\ba\label{eq:16n}
v_{l}\longrightarrow V_{0}\;\;\;\; in \;\;C^{1}_{loc}(\R^{4}_+),
\ea
On the other hand  \;$\eqref{eq:ineqvn}$\;and$\;\eqref{eq:16n}$\;imply that
\begin{equation}\label{eq:propv0n}
V_{0}(x)\leq V_{0}(0)=-\frac{1}{3}\log(k_3)\;\;\forall x\in \R^{4}_+.
\end{equation}
Moreover from\;$\eqref{eq:derbvn}$\;and\;$\eqref{eq:16n}$\;we have that \;$V_0$\;is Lipschitz.\\
On the other hand using the Green's representation formula  for\;$(P^4_g,P^3_g)$\;we obtain that for \;$x\in \R^{4}_+$\;\;fixed and for \;$R$\; big enough such that  \;$x\in B^{0}_+(R)$\;
\begin{equation}\label{eq:Gn}
u_{l}(exp_{x_{l}}(\mu_{l}x))-\bar u_{l}=\int_{M}G(exp_{x_{l}}(\mu_{l}x),y)P^{4}_{g}u_{l}(y)dV_{g}(y)+2\int_{\partial M}G(exp_{x_l}(\mu_lx),y')P^3_gu_l(y')dS_g(y').
\end{equation}
Now let us  set
\begin{equation*}
I_{l}(x)=2\int_{B^+_{x_{l}}(R\mu_{l})\cap \partial M}\left(G(exp_{x_{l}}(\mu_{l}x),y')-G(exp_{x_{l}}(0),y')\right)\bar T_{l}(y)e^{3u_{l}(y)}dS_{g}(y');
\end{equation*}
\begin{equation*}
\text{II}_{l}(x)=2\int_{\partial M \setminus( B^+_{x_{l}}(R\mu_{l})}\left(G(exp_{x_{l}}(\mu_{l}x),y')-G(exp_{x_{l}}(0),y')\right)\bar T_{l}(y')e^{3
u_{l}(y)}dS_{g}(y');
\end{equation*}
\begin{equation*}
\text{III}_{l}(x)=2\int_{\partial M}\left(G(exp_{x_{l}}(\mu_{l}x),y')-G(exp_{x_{l}}(0),y')\right) T_{l}(y)dS_{g}(y');
\end{equation*}
and
\begin{equation*}
\text{IIII}_{l}(x)=2\int_{ M}\left(G(exp_{x_{l}}(\mu_{l}x),y)-G(exp_{x_{l}}(0),y)\right) Q_{l}(y)dV_{g}(y).
\end{equation*}
Using again the same argument as in \cite{nd} (see formula (45)- formula (51)) we get
\begin{equation}\label{eq:relvn}
v_{l}(x)=I_{l}(x)+\text{II}_{l}(x)-\text{III}_{l}(x)-\text{IIII}_l(x)-\frac{1}{4}\log(3).
\end{equation}
Moreover  following the same methods as in \cite{nd}( see formula (53)-formula (62)) we obtain
\begin{equation}\label{eq:limiln}
\lim_{l}I_{l}(x)=\int_{B^{0}_+(R)\cap \partial R^4_+}\s_3\log\left(\frac{|z|}{|x-z|}\right)e^{3V_{0}(z)}dz.
\end{equation}
\begin{equation}\label{eq:limiiln}
\limsup_{l}\text{II}_{l}(x)=o_{R}(1).
\end{equation}
\begin{equation}\label{eq:limiiiln}
\text{III}_{l}(x)=o_l(1) 
\end{equation}
and
\begin{equation}\label{eq:limiiiiln}
\text{IIII}_l(x)=o_l(1).
\end{equation}
Hence from \;$\eqref{eq:16n}$,\;$\eqref{eq:relvn}$-$\eqref{eq:limiiiiln}$\; by letting \;$l$\;tends to infinity and after \;$R$\;tends to infinity, we obtain\;${V_0}_{|\R^3}$( that for simplicity we will always write by \;$V_0$)\;satisfies the following conformally invariant integral equation on \;$\R^3$
\ba\label{eq:ie1n}
V_{0}(x)=\int_{\R^3}\s_3\log\left(\frac{|z|}{|x-z|}\right)e^{3V_{0}(z)}dz-\frac{1}{3}\log(k_3).
\ea
Now since\;$V_0$\;is Lipschitz then the theory of singular integral operator gives that \;$V_0\in C^1(\R^3)$.\\
On the other hand by using the change of variable \;$y=exp_{x_l}(\mu_lx)$,\;one can check that the following holds 
\ba\label{eq:limrn}
\lim_{l\longrightarrow +\infty}\int_{B^+_{x_{l}}(R\mu_{l})\cap \partial M}\bar T_{l}e^{3u_{l}}dV_g=k_3\int_{B^+_{0}(R)\cap \partial R^4_+} e^{3V_{0}}dx;
\ea
Hence \;$\eqref{eq:limee}$\; implies that \;$e^{V_{0}}\in L^{3}(\R^{3}).$\\
Furthermore by a classification result by X. Xu, see Theorem\;$\ref{eq:xu}$\; for the solutions of \;$\eqref{eq:ie1n}$\; we derive that 
\ba\label{eq:18n}
V_{0}(x)=\log\left( \frac{2\l}{\l^{2}+|x-x_{0}|^{2}}\right)  
\ea
for some \;$\l>0 $\;\and \;$x_{0} \in \R^{3}$.\; \\ 
Moreover from  \;$V_{0}(x)\leq V_{0}(0)=-\frac{1}{3}\log(k_3)\;\;\forall x\in \R^{3}$,\;  we have that \;$\l=2k_3$\; and \;$x_{0}=0$\; namely,
$$
V_{0}(x)=\log(\frac{4\gamma_{3}}{4\gamma_{3}^{2}+|x|^{2}}).
$$
On the other hand  by letting \;$R$\;\;tends to infinity in \;$\eqref{eq:limrn}$\;we obtain
\ba\label{eq:limrln}
\lim_{R\rightarrow +\infty}\lim_{l\rightarrow +\infty}\int_{B^+_{x_{l}}(R\mu_{l})\cap \partial R^4_+}\bar T_{l}(y)e^{3u_{l}(y)}dS_{g}(y)=k_3\int_{\R^{3}}e^{3V_{0}}dx.
\ea
Moreover from a generalized Pohozaev type identity by X.Xu \cite{xu} (see Theorem 1.1) we get 
\begin{equation*}
\s_3\int_{\R^3}e^{3V_{0}(y)}dy=2,
\end{equation*}
hence using\;$\eqref{eq:limrln}$\; we derive that
$$
\lim_{R\rightarrow +\infty}\lim_{l\rightarrow +\infty}\int_{B^+_{x_{l}}(R\mu_{l})\cap\partial M}\bar T_{l}(y)e^{3u_{l}(y)}dS_{g}(y)=4\pi^2
$$\\

\noindent
Now for \;$k\geq 1$\;\;we say that \;($H_{k}$)\; holds if there exists \;$k$\;converging points\;\;$(x_{i,l})_l\subset\partial M\;\; i=1,...,k$,\;\;$k$\\sequences\;$(\mu_{i,l})\;\; i=1,...,k$\;\;of positive real numbers  converging to \;$0$\;\;such that the following hold\\ 
$
\left( A^{1}_{k}\right)$\\
$$
\hspace{-130pt}
\frac{d_{\hat g}(x_{i,l},x_{j,l})}{\mu_{i,l}}\longrightarrow +\infty \;\;\;\; i\neq j \;\;\;\;i,j =1,..,k\; and \;\;\;\;\bar T_{l}(x_{i,l})\mu_{i.l}^{3}e^{3u_{l}(x_{i,l})}=1;$$
\\
$\left( A^{2}_{k}\right)$\\
For every \;$i=1,\cdot,k$\\

$$\hspace{-130pt}x_{i,l}\rightarrow\bar x_i\in \partial M;$$
$$
v_{i,l}(x)=u_{l}(exp_{x_{i,l}}(\mu_{i,l}x))-u_{l}(x_{i,l})-\frac{1}{3}\log(k_3)\longrightarrow V_{0}(x)\;\;\;\; \text{in}\; \;\;C^{1}_{loc}(\R^4_+),\;\;\;{V_0}_{|\partial \R^4_+}:=\log(\frac{4\gamma_{3}}{4\gamma_{3}^{2}+|x|^{2}})
$$
and
$$
\hspace{-110pt} 
\lim_{R\rightarrow +\infty}\lim_{l\rightarrow +\infty}\int_{B^+_{x_{i,l}}(R\mu_{i,l})\cap\partial M}\bar T_{l}(y)e^{3u_{l}(y)}=4\pi^2
$$

\noindent
Clearly, by the above arguments\;$(H_{1})$\;holds. We let now \;$k\geq 1$\; and assume that \;$(H_{k})$\; holds. We also assume that
\ba \label{eq:lsR}
\sup_{\partial M}R_{k,l}(x)^{3}e^{3u_{l}(x)}\longrightarrow+\infty\;\;\; as \;\;\;l\longrightarrow +\infty,
\ea
where$$R_{k,l}(x)=\min_{i=1;..;k}d_{g}(x_{i,l},x).$$
Now using the same argument as in  \cite{dr},\cite{nd} and the arguments which have rule out the possibility of interior blow up above that also apply for local maxima, one can see easily that $(H_{k+1})$.
Hence  since \;$\left( A^{1}_{k}\right) $\; and \;$\left( A^{2}_{k}\right) $\; of \;$H_{k}$\; imply that 
\begin{equation*}
\int_{\partial M}\bar T_l(y)e^{3u_{l}(y)}dS_g(y)\geq k4\pi^2+o_l(1).
\end{equation*}
Thus \;$\eqref{eq:limee}$\;imply that there exists a maximal \;$k$\;,\;$1\leq k\leq \frac{1}{4\pi^2}\left( \int_{M}Q_{0}(y)dV_g(y)+\int_{\partial M}T_0(y')dS_g(y')\right)$\;, such that \;$(H_{k})$\;holds. Arriving to this maximal \;$k$, we get that \;$\eqref{eq:lsR}$\; cannot hold. Hence setting \;$N=k$\; the proof of Step 1 is done.\\

{\bf Step 2}\\
\noindent
There exists a constant \;$C>0$\; such that
\ba\label{eq:grad}
R_{l}(x)|\n_{g} u_l|_{g}(x)\leq C\;\;\;\;\;\;\;\;\;\;\;\;\;\forall x\in M\;\;and\;\;\forall l\in N;
\;\;\;\forall\; x\in \partial M
\ea
where 
$$R_{l}(x)=\min_{i=1,..,N}d_{g}(x_{i,l},x);
$$
and the \;$x_{i,l}$'s are as in Step 1.\\
\vspace{10pt}

\noindent
{\em Proof of Step 2}\\
\vspace{4pt}
First of all using the Green representation formula for \;$(P^4_g,P^3_g)$\;see Lemma\;$\ref{eq:greenrep}$\; we obtain
\begin{equation*}
u_l(x)-\bar u_l=\int_{M}G(x,y)P^4_gu_l(y)dV_g(y)+2\int_{\partial M}G(x,y')P^3_gu_l(y')dS_g(y').
\end{equation*}
Now using the BVP\;$\eqref{eq:bvps}$\;we get
\begin{equation}\label{eq:todif}
\begin{split}
u_l(x)-\bar u_l=-2\int_{M}G(x,y)Q_ldV_g(y)-2\int_{\partial M}G(x,y')T_l(y')u_l(y')dS_g(y')\\+2\int_{\partial M}G(x,y)\bar T_l(y')e^{3u_l(y')}dS_g(y').
\end{split}
\end{equation}
Thus differentiating with respect to \;$x$\;$\eqref{eq:todif}$\; and using the fact that \;$Q_l\rightarrow Q_0$,\;$\bar Q_l\rightarrow \bar Q_0$\;and \;$T_l\rightarrow T_0$\; in \;$C^1$, we have that for \;$x_l\in \partial M$
\noindent
$$
|\n_{ g} u_l(x_l)|_{g}=O\left( \int_{\partial M }\frac{1}{d_g(x_l,y)} e^{3u_l(y)}dS_g(y)\right) +O(1).
$$
Hence at this stage following the same argument as in the proof of Theorem 1.3, Step 2 in \cite{nd},  we obtain 
\begin{equation*}
\int_{\partial M}\frac{1}{(d_{g}(x_{l},y))}e^{3u_l(y)}dV_{g}(y)=O\left( \frac{1}{R_l(x_l)}\right);
\end{equation*}
hence  since \;$x_l$\;is arbitrary, then the proof of Step 2 is complete.
\vspace{10pt}

\noindent
{\bf Step 3}\\
\noindent
Set 
\begin{equation*}
R_{i,l}=\min_{i\neq j}d_{g}(x_{i,l},x_{j,l});
\end{equation*}
we have that 
\vspace{2pt}
 
\noindent
1) There exists a constant\;\;$ C>0$\; such that \;\;$\forall \;r\in (0,R_{i,l}]\;\;\forall \; s\in
(\frac{r}{4},r]$\\
\ba\label{eq:harn2}
|u_{l}\left( exp_{x_{i,l}}(rx)\right) -u_{l}\left( exp_{x_{i,l}}(sy)\right)|\leq C \;\;\;\;
for\;\; all \;\;x,y \in \partial \R^4_+ \;\;such \;\;that \;\;|x|,\;|y|\leq \frac{3}{2}.
\ea
2)\;If \;\;$d_{i,l}$\;is such that \;$0<d_{i,l}\leq \frac{R_{i,l}}{2}$\; and\;$\frac{d_{i,l}}{\mu_{i,l}}\longrightarrow+\infty$\; then we have that\\
if
\ba\label{eq:d_{i,l}2}
\int _{B^+_{x_{i,l}}(d_{i,l})\cap \partial M}\bar T_{l}(y)e^{3u_l(y)}dS_{g}(y)=4\pi^2+o_{l}(1);
\ea
then
$$
\int_{B^+_{x_{i,l}}(2d_{i,l})\cap \partial M}\bar T_{l}(y)e^{3u_{l}(y)}ds_{g}(y)=4\pi^2+o_{l}(1).
$$ 
3) Let \;$R$\;be large and fixed. If \;$d_{i,l}>0$\;is such that \;$d_{i,l}\longrightarrow 0$,\;$\frac{d_{i,l}}{\mu_{i,l}}\longrightarrow+\infty$, and \;$d_{i,l}<\frac{R_{i,l}}{4R}$\\
then 
if 
$$\int_{B^+_{x_{i,l}}(\frac{d_{i,l}}{2R})\cap \partial M}\bar Q_{l}(y)e^{3u_{l}(y)}dS_{g}(y)=4\pi^2+o_{l}(1);$$\;
then by setting
$$
\tilde u_{l}(x)=u_{l}(exp_{x_{i,l}}(d_{i,l}x));\;\;\;\;\;\;\;\;\;x\in  A^+_{2R};
$$
where \;$A^+_{2R}=(B^{0}_+(2R)\setminus B^{0}_+(\frac{1}{2R}))\cap\partial \R^4_+$,\;we have that, 
$$
||d_{i,l}^{4}e^{3\tilde u_{l}}||_{C^{\alpha}(A^+_{R})}\rightarrow 0\;\;as\;\;l\rightarrow +\infty;
$$
for some \;$\alpha \in (0,1)$\;where \;$A^+_{R}=(B^{0}_+(R)\setminus B^{0}_+(\frac{1}{R}))\cap\partial \R^4_+$.\\

\noindent
{\em Proof of Step 3}\\
\noindent
We have that property 1 follows immediately from Step 2 and the definition of \;$R_{i,l}$.  In fact we can join \;$rx$\; to\;$sy$\; by a curve whose length is bounded by a constant proportional to \;$r$.\\
Now let us show  point 2. Thanks to\;$\frac{d_{i,l}}{\mu_{i,l}}\longrightarrow+\infty$,\;point c) of Step 1 and \;$\eqref{eq:d_{i,l}2}$\;we have that 
\ba\label{eq:ann2}
\int_{B^+_{x_{i,l}}(d_{i,l})\cap\partial M\setminus B^+_{x_{i,l}}(\frac{d_{i,l}}{2})\cap \partial M}e^{3u_l(y)}dS_{g}(y)=o_{l}(1).
\ea 
Thus   using \;$\eqref{eq:harn2}$, with \;$s=\frac{r}{2}$\;and\;$r=2d_{i,l}$\;we get
$$
\int_{B^+_{x_{i,l}}(2d_{i,l})\cap \partial M\setminus B^+_{x_{i,l}}(d_{i,l})\cap \partial M}e^{3u_{l}(y)}ds_g(y)\leq C\int_{B^+_{x_{i,l}}(d_{i,l})\cap \partial M\setminus B^+_{x_{i,l}}(\frac{d_{i,l}}{2})\cap \partial M}e^{3u_l(y)}dS_{g}(y);
$$
Hence we arrive 
$$
\int_{B^+_{x_{i,l}}(2d_{i,l})\cap \partial M\setminus B^+_{x_{i,l}}(d_{i,l})\cap \partial M}e^{3u_{l}(y)}dS_g(y)=o_{l}(1).
$$
So the proof of point 2 is done.
On the other hand by following in a straightforward way the proof of point 3 in Step 3 of Theorem 1.3 in \cite{nd} one gets easily point 3. Hence the proof of Step 3 is complete.
\vspace{10pt}

\noindent
{\bf Step 4}\\
\noindent
There exists a positive constant \;$C$\; independent of \;$l$\; and\;$i$\; such that
\begin{equation*}
\int_{B^+_{x_{i,l}}(\frac{R_{i,l}}{C})\cap \partial M}\bar T_{l}(y)e^{3u_{l}(y)}dS_g(y)=4\pi^2+o_{l}(1).
\end{equation*}

\vspace{6pt}

\noindent
{\em Proof of Step 4}\\
\noindent
The proof is an adaptation of the arguments in Step 4 (\cite{nd})
\vspace{20pt}

\noindent
{\bf Step 5 :Proof of Theorem\;$\ref{eq:comp}$}\\
\vspace{6pt}

\noindent
Following the same argument as in Step 5(\cite{nd}) we have
\begin{equation*}
\int_{\partial M \setminus (\cup_{i=1}^{i=N}B^+_{x_{i,l}}(\frac{R_{i,l}}{C})\cap\partial M)}e^{3u_{l}(y)}dS_{g}(y)=o_{l}(1).
\end{equation*}

So since \;$B^+_{x_{i,l}}(\frac{R_{i,l}}{C})\cap\partial M$\; are disjoint then the Step 4 implies that, 
$$
\int_{\partial M}\bar T_{l}(y)e^{3u_{l}(y)}dS_{g}(y)=4N\pi^2+o_{l}(1),
$$
hence \;$\eqref{eq:limee}$\;implies that 
$$
\int_{ M} Q_{0}(y)dV_{g}(y)+\int_{\partial M}T_0(y')dS_g(y')=4N\pi^2.
$$
ending the proof of Theorem\;$\ref{eq:comp}$.
\end{pfn}

\section{Proof of Theorem\;$\ref{eq:theorem3}$}
This section deals with the proof of Theorem\;$\ref{eq:theorem3}$. It is divided into four Subsections. The first one is concerned with an improvement of the Moser-Trudinger type inequality (see Proposition\;$\ref{eq:mos-trub}$) and its corollaries. The second one is about the existence of a non-trivial global projection from some negative sublevels of \;$II$\;onto \;$\partial M_k$\;(for the definition see Section 2 formula\;$\ref{eq:base}$). The third one deals with the construction of a map from \;$\partial M_k$\;into suitable negative sublevels of \;$II$. The last one describes the min-max scheme.
\subsection{Improved Moser-Trudinger inequality}
In this Subsection we give an improvement of the Moser-Trudinger type inequality, see Proposition\;$\ref{eq:mos-trub}$.  Afterwards, we state a Lemma which gives some sufficient conditions for the improvement to hold (see \;$\eqref{eq:rat}$). By these results, we derive that, for  \;$u\in H_{\frac{\partial }{\partial n}}$\; such that\;$II(u)$\;attains large negative values, \;$e^{3u}$\;can concentrate at most at \;$k$\;points of \;$\partial M$. (see Lemma\;$\ref{eq:poc}$). Finally from these results, we derive a corollary which gives the distance of \;$e^{3u}$\;(for some functions \;$u$\; suitably normalized) from \;$\partial M_k$.
\vspace{6pt}

\noindent
As said in the introduction of the Subsection, we start by the following Lemma giving an improvement of the Moser-Trudinger type inequality (Proposition\;$\ref{eq:mos-trub}$). Its proof is a trivial adaptation of the arguments of Lemma 2.2in \cite{dm}.
\begin{lem}\label{eq:impro}
For a fixed \;$l\in\N$,  let \;$S_{1}\cdots S_{l+1}$, be subsets of\;$\partial M$\;satisfying,\;$dist(S_i,S_j)\geq \d_0$\;for \;$i\neq j$, let \;$\gamma_0\in (0,\frac{1}{l+l})$.\\
Then, for any\;$\bar \epsilon>0$,\;there exists a constant\;$C=C(\bar \epsilon,\d_0,\gamma_0,l,M,)$\;such that  the following hods\\
1)
\begin{equation*}
\log\int_{\partial M}e^{3(u-\bar u_{\partial M})}\leq C+\frac{3}{16\pi^2}(\frac{1}{l+1-\bar \epsilon})\left<P^{4,3}_{g}u,u\right>_{L^2(M)};
\end{equation*}
for all the functions \;$u\in H_{\frac{\partial }{\partial n}}$\;satisfying
\begin{equation}\label{eq:rat}
\frac{\int_{S_i}e^{3u}dSg}{\int_{\partial M}e^{3u}dSg}\geq \gamma_0,\;\;i\in\lbrace1,..,l+1\rbrace.
\end{equation}

\end{lem}
\vspace{6pt}

\noindent
In the next Lemma we show a criterion which implies the situation described in the first condition in \;$\eqref{eq:rat}$. The result is proven in\;$\cite{dm}$\;Lemma 2.3.

\begin{lem}\label{eq:er}
Let \;$l$\;be a given positive integer, and suppose that \;$\epsilon$\;and \;$r$ are positive numbers. Suppose  that for a non-negative function\;$f \in L^1(\partial M)$\;with
$\|f\|_{L^1(\partial M)} = 1$\;there holds
\begin{equation}\nonumber
  \int_{\cup_{i=1}^\ell B^{\partial M}_r(p_i)} f dS_g < 1 - \e \qquad \qquad \hbox{ for every\; $\ell$-tuples }
  p_1, \dots, p_\ell \in \partial M
\end{equation}
Then there exist $\ov{\e} > 0$ and $\ov{r} > 0$, depending only on
$\e, r, \ell$ and $\partial M$ (but not on $f$), and $\ell + 1$ points
$\ov{p}_1, \dots, \ov{p}_{\ell+1} \in \partial M$ (which depend on $f$)
satisfying
$$
  \int_{B^{\partial M}_{\ov{r}}(\ov{p}_1)} f dS_g > \ov{\e}, \; \dots, \;
  \int_{B^{\partial M}_{\ov{r}}(\ov{p}_{\ell+1})} f dS_g > \ov{\e}; \qquad \qquad
  B^{\partial M}_{2 \ov{r}}(\ov{p}_i) \cap B^{\partial M}_{2 \ov{r}}(\ov{p}_j) = \emptyset
  \hbox{ for } i \neq j.
$$
\end{lem}
\vspace{4pt}

\indent

\indent
An interesting consequence of Lemma\;$\ref{eq:impro}$\;is the following one. It characterize some functions in \;$H_{\frac{\partial}{\partial N}}$\;for which the value of \;$II$\;is large negative. 
\begin{lem}\label{eq:poc}
Under the assumptions of Theorem\;$\ref{eq:theorem3}$, and for \;$k\geq 1$\;given by \;$\eqref{eq:range}$, the following property holds. For any \;$\epsilon>0$\;and any \;$r>0$\;there exists large positive\;$L=L(\epsilon,r)$\;such that for any \;$u\in H_{\frac{\partial }{\partial n}}$\;\;with\;\;$II(u)\leq -L,\;\;\int_{\partial M}e^{3u}dS_{g}=1$\;\;there exists\;\;$k$\;\;points \;\;$p_{1,u},\dots,p_{k,u}\in \partial M$\;\;such that
\begin{equation}
\int_{\partial M\setminus \cup_{i=1}^{k}B^{\partial M}_{p_{i,u}}(r)}e^{3u}dS_{g}<\epsilon
\end{equation}
\end{lem}
\begin{pf}
Suppose that by contradiction the statement is not true. Then there exists \;$\epsilon>0$,\;$r>0$, and a sequence \;$(u_n)\in H_{\partial n}$\; such that \;$\int_{\partial M}e^{3u_n}dS_g=1$, \;$II(u_n)\rightarrow -\infty$\;as\;$n\rightarrow+\infty$ \;and such that\\
for any \;$ k$\;tuples of points \;$p_{1},\dots,p_{k}\in \partial M$\;,we have 
\begin{equation}\label{eq:epsilones}
\int_{( \cup_{i=1}^{k}B^{\partial M}_{p_{i,u}}(r))}e^{3u}dS_{g}<1-\epsilon;
\end{equation}
Now  applying Lemma\;$\ref{eq:er}$\; with \;$f=e^{3u_n}$, and after Lemma\;$\ref{eq:impro}$\;with \;$\delta_0=2\bar r$, \;$S_i=B^{\partial M}_{\bar p_i}(\bar r)$, \;and \;$\gamma_0=\bar \epsilon$\; where \;$\bar\epsilon$\;,\;$\bar r$,\;$\bar p_i$\;\;are given as in Lemma\;$\ref{eq:er}$, we have for every \;$\tilde \epsilon>0$\; there exists\;$C$\;depending on \;$\epsilon$,\;$r$, and \;$\tilde \epsilon $\;such that 
\begin{equation*} 
\begin{split}
II(u_n)\geq \left<P^{4,3}_gu_n,u_n\right>+4\int_{M}Q_gu_ndV_g+4\int_{\partial M}T_gu_ndS_g-\frac{4}{3}\kappa_{(P^4,P^3)}\frac{3}{16\pi^2(k+1-\tilde \epsilon)}\left<P^{4,3}_gu_n,u_n\right>\\-C\kappa_{(P^4,P^3)}-4\kappa_{(P^4,P^3)}\ov{u_n}_{\partial M}
\end{split}
\end{equation*}
where  $C$\;is independent of \;$n$.
Using elementary simplifications, the above inequality becomes 
\begin{equation*} 
\begin{split}
II(u_n)\geq \left<P^{4,3}_gu_n,u_n\right>+4\int_{M}Q_gu_ndV_g+4\int_{\partial M}T_gu_ndS_g-\frac{\kappa_{P^4,P^3}}{4\pi^2(k+1-\tilde \epsilon)}\left<P^{4,3}_gu_n,u_n\right>\\-C\kappa_{P^4,P^3}-4\kappa_{P^4,P^3}\ov{u_n}_{\partial M}.
\end{split}
\end{equation*}
So, since \;$\kappa_{P^4,P^3}<(k+1)4\pi^2$, by choosing  \;$\tilde \epsilon$\;small we get 
\begin{equation*}
II(u_n)\geq \beta\left<P^{4,3}_gu_n,u_n\right>-4C\left<P^{4,3}_gu_n,u_n\right>^{\frac{1}{2}}-C\kappa_{P^4,P^3};
\end{equation*}
thanks to  H\"{o}lder inequality, to Sobolev embedding,  to trace Sobolev embedding and to the fact that \;$Ker P^{4,3}_{g_0}\simeq \R$ (where \;$\beta=1-\frac{\kappa_{P^4,P^3}}{4\pi^2(k+1-\tilde \epsilon)}>0$). Thus we arrive to 
\begin{equation*}
II(u_n)\geq -C.
\end{equation*}
So we reach a contradiction. Hence the Lemma is proved.
\end{pf}

\vspace{4pt}

\noindent

\noindent 
Next we give a Lemma which is a direct consequence of the previous one. It gives the distance of the functions $e^{3u}$,
 from $\partial M_k$ \; for \;$u$\;belonging to low energy levels of \;$II$\; such that \;$\int_{\partial M}e^{3u}dS_g=1$. Its proof is the same as the one of corollary in \cite{dm}.

\begin{cor}\label{eq:impmt}
Let $\ov{\e}$ be a (small) arbitrary positive number and \;$k$\;be given as in\;$\eqref{eq:range}$. Then there
exists $L > 0$ such that, if \;$\text{II}(u) \leq - L$ \;and \;
$\int_{\partial M} e^{3u} dS_g = 1$, then  we have that\;$d(e^{3u}, \partial M_k) \leq 
\ov{\e}$.
\end{cor}

\subsection{ Mapping sublevels of \;$II$\;into \;$(M_{\partial})_k$}
In this short Subsection we show that one can map in a non trivial way some appropriate low energy sublevels of the Euler-Lagrange functional \;$II$\;into \;$\partial M_k$.\\
First of all arguing as in Proposition 3.1 in \;$\cite{dm}$, we have the following Lemma.

\begin{lem}\label{eq:projdm} Let $m$ be a positive integer, and for
$\e > 0$ let $\mathcal{D}_{\e,m}$ be as in \;$\eqref{eq:dm}$. Then there exists
$\e_m > 0$, depending on $m$ and $\partial M$ such that, for $\e \leq
\e_k$ there exists a continuous map $\Pi_m : \mathcal{D}_{\e,m} \to
\partial M_m$.
\end{lem}
Using the above Lemma we have the following non-trivial continuous global projection form low energy sublevels of \;$II$\;into \;$\partial M_k$.
\begin{pro}\label{eq:pro1}
For \;$k\geq 1$\;given as in \;$\eqref{eq:range}$,\;there exists a large \;$L>0$ and a continuous map \;$\Psi $\;from the sublevel\;\;$\lbrace u: II(u)<-L,\;\int_{\partial M}e^{3u}dS_g=1\rbrace$\;into \;$\partial M_{k}$\;which is topologically non-trivial. 
\end{pro}

\noindent
By the non-contractibility of \;$\partial M_k$, the non-triviality of the map is apparent from \;{\em b)}\; of Proposition\;$\ref{eq:pro2}$\;below.

\begin{pf}
We fix $\e_{k}$ so small that Lemma\;$\ref{eq:projdm}$\; applies
with $m = k$. Then we apply Corollary \ref{eq:impmt} with $\ov{\e} =
\e_k$. We let $L$ be the corresponding large number, so that if
$II(u) \leq - L$\; and \;$\int_{\partial M}e^{3u}dS_g=1$, then $d(e^{3u}, \partial M_{k}) < \e_k$. Hence
for these ranges of $u$ , since the map $u \mapsto e^{3u}$
is continuous from $H^1(M)$ into $L^1(\partial M)$, then the projections
$\Pi_{k}$ from $H^1(\Sig)$ onto $\partial M_k$ is well defined and
continuous. .
\end{pf}

\subsection{Mapping \;
$\partial M_k$\;into sublevels of\;$II$} 
In this Subsection we will define some test functions depending on a real parameter\;$\l$\;and give estimate of the quadratic part of the functional \;$II$\;on those functions as \;$\l$\; tends to infinity. And as a corollary we define  a continuous map from \;$\partial M_{k}$\; into large negative sublevels of \;$II$.\\
For\;$\d>0$\;small, consider a smooth non-decreasing cut-off function \;$\chi_{\delta}:\R_{+}\rightarrow \R$\;\;satisfying the following properties (see \;$\cite{dm}$): 
\begin{equation*}
\left\{
\begin{array}{ll}
\chi_{\d}(t)=t, \;\;\;\;for\;\;t\in[0,\d];\\
\chi_{\d}(t)=2\d, \;\;\;\;for\;\;t\geq 2\d;\\
\chi_{\d}(t)\in [\d,2\d ],\;\;\;\;for\;\;t\in[\d,2\d].
\end{array}
\right.
\end{equation*}
Then, given \;$\sigma=\in \partial M_k$, $\sigma=\sum_{i=1}^{k}t_i\d_{x_i}$\; and \;$\l>0$, we define the function ;$\varphi_{\l,\sigma}:M\rightarrow\ \R$\;as follows
\begin{equation}\label{eq:vsl}
\varphi_{\l,\sigma}(y)=\frac{1}{3}\log\left[ \sum_{i=1}^{k }t_i\left(\frac{2\l}{1+\l^2\chi_{\d}^{2}(d_{i}(y))}\right)^3\right];
\end{equation}
where we have set
\begin{equation}\nonumber
d_i(y)=d_g(y,x_i),\;\;\;\;\;\;\;x_i\in \partial M,y\in M,;
\end{equation}

with\;$d_g(\cdot,\cdot)$\;denoting the Riemannian distance on \;$M$. \\
 Now we state a Lemma giving an estimate (uniform in \;$\s\in \partial M_k$) of the quadratic part \;$\left<P^{4,3}_{g}\varphi_{\l,\sigma},\varphi_{\l,\sigma}\right>$\;of the Euler functional \;$II$\;as \;$\l\rightarrow +\infty$. Its proof is a straightforward adaptation of the arguments in Lemma 4.5 in \cite{nd}. 
\begin{lem}\label{eq:qap}
Suppose \;$\varphi_{\l,\sigma}$\;as in \;$\eqref{eq:vsl}$\; and let \;$\epsilon>0$\;small enough. Then as \;$\l\rightarrow +\infty$\;one has 
\begin{equation}\label{eq:esmai}
\left<P^{4,3}_{g}\varphi_{\l,\sigma},\varphi_{\l,\sigma}\right>\leq (16\pi^2k+\epsilon+o_{\delta}(1))\log \l +C_{\epsilon,\delta}
\end{equation}
\end{lem}
\vspace{6pt}

\noindent
Next we state a lemma giving estimates of the remainder part of the functional \;$II$\;along \;$\varphi_{\s,\l}$. The proof is the same as the one of formulas (40) and (41) in the proof of Lemma 4.3 in \;$\cite{dm}$.
\begin{lem}\label{eq:lep}
Soppose \;$\varphi_{\s,\l}$\;as in \;$\eqref{eq:vsl}$.\;Then as \;$\l\rightarrow+\infty$\;one has
\begin{equation*}
\int_{M}Q_g\varphi_{\s,\l}dV_g=-\kappa_{P^{4}_g}\log \l+O(\d^4\log \l)+O(\log \d)+O(1);
\end{equation*}
\begin{equation*}
 \int_{\partial M}T_g\varphi_{\s,\l}dV_g=-\kappa_{P^{3}_g}\log \l+O(\d^3\log \l)+O(\log \d)+O(1);
\end{equation*}
and
\begin{equation*}
\log \int_{\partial M}e^{3\varphi_{\s,\l}}=O(1).
\end{equation*}
\end{lem}
Now for \;$\l>0$\;we define the map \;$\Phi_{\l}:\partial M_k\rightarrow H_{\frac{\partial }{\partial n}}$\;by the following formula
\begin{equation*}
\forall \;\;\s\in \partial M_k\;\;\;\Phi_{\l}(\s)=\varphi_{\s,\l}.
\end{equation*}
We have the following Lemma which is a trivial application of Lemmas\;$\ref{eq:qap}$\;and\;$\ref{eq:lep}$.
\begin{lem}\label{eq:proj}
For \;$k\geq 1$\;(given as in \;$\eqref{eq:range}$\;), given any \;$L>0$\;large enough, there exists a small \;$\d$\;and a large \;$\bar \l$\;such that \;$II(\Phi_{\bar \l}(\s))\leq -L$\;for every \;$\s\in \partial M_k$.
\end{lem}
Next we state a proposition giving the existence of the projection from \;$\partial M_k$\;into large negative sublevels of \;$II$, and the non-triviality of the map \;$\Psi$\;of the proposition\;$\eqref{eq:pro1}$.
\begin{pro}\label{eq:pro2}
Let \;$\Psi$\;be the map defined in proposition\;$\ref{eq:pro1}$\;. Then assuming \;$k\geq 1$\;(given as in \;$\eqref{eq:range}$),\;for every \;$L>0$\;sufficiently large (such that proposition\;$\ref{eq:pro1}$\;applies), there exists a map
\begin{equation*}
\Phi_{\bar \l}:\partial M_k\longrightarrow H_{\frac{\partial }{\partial n}}\;\;   
\end{equation*}
with the following properties\\
a)
\begin{equation*}
II(\Phi_{\bar \l}(z))\leq -L\;\;\text{for any}\;\;z\in \partial M_{k}; 
\end{equation*}
b)\\
$\Psi \circ\Phi_{\bar \l}$\;is homotopic to the identity on \;$\partial M_{k}$. 
\end{pro}
\vspace{6pt}

\noindent
\begin{pf}
The statement {\em(a)} follows from Lemma \;$\ref{eq:proj}$. To prove {\em(b)} it is sufficient to consider the family of maps \;$T_{\l}:\partial M_k\rightarrow \partial M_k$\;defined by 
\begin{equation}\nonumber
T_{\l}(\s)=\Psi(\Phi_{\l}(\s)),\;\;\;\;\s\in \partial M_k
\end{equation}
We recall that when \;$\l$\;is sufficiently large, then this composition is well defined. Therefore , since \;$\frac{e^{3\varphi_{\s,\l}}}{\int_{\partial M}e^{3\varphi_{\s,\l}}dS_g}\rightharpoonup \s$\;in the weak sens of distributions, letting \;$\l\rightarrow +\infty$\;we obtain an homotopy between \;$\Psi\circ\Phi$\;and\;$\text{Id}_{\partial M_k}$. This concludes the proof.
\end{pf}
\vspace{20pt}

\subsection{Min-max scheme}
In this Subsection, we describe the min-max scheme based on the set \;$\partial M_{k}$\;in order to prove Theorem\;$\ref{eq:theorem3}$. As anticipated in the introduction, we define a modified functional \;$II_{\rho}$\; for which we can prove existence of solutions in a dense set of the values of \;$\rho$. Following a idea of Struwe ( see \;$\cite{str}$), this is done by proving the a.e differentiability of the map \;$\rho \rightarrow \ov{ II}_{\rho}$\;( where \;$\ov{ II}_{\rho}$\;is the minimax value for the functional \;$II_{\rho}$).\\
\vspace{6pt}

\noindent
We now introduce the minimax scheme which provides existence of solutions for (8). Let \;$\widehat {\partial M_k}$\;denote the (contractible) cone over\;$\partial M_{k}$, which can be represented as \;$\widehat {\partial M_{k}}=(\partial M_{k}\times[0,1])$\;with\;$\partial M_{k}\times\text{{0}}$\;collapsed to a single point. First let\;$L$\;be so large that Proposition\;$\ref{eq:pro1}$\; applies with \;$\frac{L}{4}$, and then let \;$\bar \l$\;be so large that Proposition\;$\ref{eq:pro2}$\;applies for this value of \;$L$. Fixing\;$\bar \l $, we define the following class.
\begin{equation}\label{eq:class}
II_{\bar \l}=\lbrace \pi:\widehat {\partial M_{k}}\rightarrow H_{\frac{\partial}{\partial n}}:\pi\;\text{is continuous and }\;\pi(\cdot\times\text{{1}})=\Phi_{\bar \l}(\cdot)\rbrace.
\end{equation}
We then have the following properties.
\begin{lem}\nonumber
The set\;$II_{\bar \l}$\;is non-empty and moreover, letting
\begin{equation}\nonumber
\ov {II}_{\bar \l}=\inf_{\pi\in II_{\bar \l}}\sup_{m\in \widehat {\partial M_{k,}}}II(\pi(m)),\;\;\;\;\;\text{there holds}\;\;\;\ov {II}_{\bar \l}>-\frac{L}{2}.
\end{equation}
\end{lem}
\begin{pf}
The proof is the same as the one of Lemma 5.1 in \;$\cite{dm}$. But we will repeat it for the reader's convenience.\\
To prove that \;$\ov {II}_{\bar \l}$\;is non-empty, we just notice that the following map 
\begin{equation*}
\bar \pi(\cdot,t)=t\Phi_{\bar \l}(\cdot)
\end{equation*}
belongs to \;$II_{\bar \l}$. Now to prove that \;$\ov {II}_{\bar \l}>-\frac{L}{2}$,\;let us argue by contradiction. Suppose that \;$\ov {II}_{\bar \l}\leq-\frac{L}{2}$: then there exists a map \;$\pi\in {II}_{\bar \l}$\; such that \;$\sup_{m\in \widehat {\partial M_{k}}}II(\pi(m))\leq-\frac{3}{8}L$. Hence since Proposition\;$\ref{eq:pro1}$\; applies with \;$\frac{L}{4}$, writing\;$ m=(z,t)$\;with \;$z\in \partial M_{k}$\;we have that the map 
\begin{equation}\nonumber
t\rightarrow \Psi\circ\pi(\cdot,t)
\end{equation}
is an homotopy in \;$\partial M_{k}$\;between \;$\Psi\circ\Phi_{\bar \l}$\;and a constant map. But this is impossible since \;$\partial M
_{k}$\;is non-contractible and \;$\Psi\circ\Phi_{\bar \l}$\;is homotopic to the identity by Proposition\;$\ref{eq:pro2}$.\\
\end{pf}
\vspace{6pt}

\noindent
Next we introduce a variant of the above minimax scheme, following\;\;$\cite{dm}$\;$\cite{str}$\;and\cite{nd}. For \;$\rho$\;in a small neighborhood of \;$1$, \;$[1-\rho_{0}, 1+\rho_{0}]$, we define the modified functional \;$II_{\rho}:H_{\frac{\partial}{\partial n}}\rightarrow \R$
\begin{equation}\label{eq:pertf}
II_{\rho}(u)=\left<P^{4,3}_{g}u,u\right>+4\rho\int_{M}Q_{g}udV_{g}+4\rho\int_{\partial M}T_gudS_g-\frac{4}{3}\rho\kappa_{(P^4,P^3)}\log\int_{\partial M}e^{3u}dS_{g};\;\;\;u\in H_{\frac{\partial}{\partial n}}.
\end{equation}
Following the estimates of the previous section, one easily checks that the above minimax scheme applies uniformly for \;$\rho\in[1-\rho_{0}, 1+\rho_{0}]$\;and for \;$\bar \l$\;sufficiently large. More precisely, given any large number \;$L>0$, there exist \;$\bar \l$\;sufficiently large and \;$\rho_{0}$\;sufficiently small such that 
\begin{equation}\label{eq:espef}
 \sup_{\pi\in II_{\bar \l}}\sup_{m\in \partial\widehat {\partial M_{k}}}II(\pi(m))<-2L;\;\;\;\ov {II}_{\rho}\inf_{\pi\in II_{\bar \l}}\sup_{m\in \widehat {\partial M_{k}}}II_{\rho}(\pi(m))>-\frac{L}{2};\;\;\;\;\rho\in[1-\rho_{0},1+\rho_{0}],
\end{equation}
where\;$II_{\bar \l}$\;is defined as in \;$\eqref{eq:class}$. Moreover, using for example the test map, one shows that for \;$\rho_{0}$ \;sufficiently small there exists a large constant \;$\bar L$\; such that 
\begin{equation}\label{eq:esrb}
\ov {II}_{\rho}\leq \bar L,\;\,\;\;\;\;\text{for every}\;\rho\in[1-\rho_{0},1+\rho_{0}].
\end{equation}
\vspace{6pt}

\noindent
We have the following result regarding the dependence in \;$\rho$\;of the minimax value \;$\ov{II}_{\rho}$. 
\begin{lem}\label{eq:noinc}
Let \;$\bar \l$\; and \;$\rho_{0}$\;such that \;$\eqref{eq:espef}$\;holds. Then the function 
\begin{equation}\nonumber
\rho\rightarrow \frac{\ov {II}_{\rho}}{\rho}\;\;\;\;\;\;\text{is non-increasing in}\;\;\;[1-\rho_{0},1+1-\rho_{0}]
\end{equation}
\end{lem}
\begin{pf}
For\;$\rho \geq \rho^{'}$, there holds
\begin{equation*}
\frac{II_{\rho}(u)}{\rho}-\frac{ II_{\rho^{'}}(u)}{\rho^{'}}=\left(\frac{1}{\rho}-\frac{1}{\rho^{'}}\right)\left<P^{4,3}_gu ,u\right>
\end{equation*}
Therefore it follows easily that also
\begin{equation*}
\frac{\ov{II}_{\rho}}{\rho}-\frac{ \ov{II}_{\rho^{'}}}{\rho^{'}}\leq 0,
\end{equation*}
hence the Lemma is proved.
\end{pf}
\vspace{6pt}

\noindent
From this Lemma it follows that the function \;$\rho\rightarrow\frac{\ov {II}_{\rho}}{\rho}$\;is a.e. differentiable in \;$[1-\rh_0,1+\rho_0]$, and we obtain the following corollary.
\begin{cor}\label{eq:cra}
Let \;$\bar \l$\;and\;$\rho_{0}$\;be as in Lemma\;$\ref {eq:noinc}$, and let\;$\Lambda\subset    [1-\rho_{0},1+\rho_{0}]$\;be the (dense) set of \;$\rho$\;\;for which the function\;\;$\frac{\ov {II}_{\rho}}{\rho}$\;\;is differentiable. Then for \;$\rho\in \Lambda$\; the functional \;$II_{\rho}$\;possesses a bounded Palais-Smale sequence \;$(u_{l})_{l}$\;at level \;$\ov {II}_{\rho}$.
\end{cor}
\vspace{4pt}

\noindent
\begin{pf}
The existence of Palais-Smale sequence \;$(u_{l})_{l}$\; at level \;$\ov {II}_{\rho}$\;follows from \;$\eqref{eq:espef}$\; and the bounded is proved exactly as in \cite{djlw}, Lemma 3.2.
\end{pf}
\vspace{6pt}

\noindent
Next we state a Proposition saying that bounded Palais-Smale sequence of \;$II_{\rho}$ converges weakly (up to a subsequence) to a solution of the perturbed problem. The proof is the same as the one of Proposition 5.5 in \;$\cite{dm}$.
\begin{pro}\label{eq:pss}
Suppose \;$(u_{l})_{l} \subset H_{\frac{\partial}{\partial n}}$\;is a sequence for which 
\begin{equation}\nonumber
II_{\rho}(u_{l})\rightarrow c\in \R;\;\;\;\;\;II^{'}_{\rho}[u_{l}]\rightarrow 0;\;\;\; \int_{\partial M}e^{3u_l}dS_g=1\;\;\;\|u_l\|_{H^2(M)}\leq C.
\end{equation}
Then \;$(u_{l})$\;has a weak limit \;$u$\;(up to a subsequence) which satisfies the following equation:
\begin{equation}\nonumber
\left\{
\begin{split}
P^4_gu+2\rho Q_g&=0\;\;&\text{in}\;\;M;\\
P^3_gu+\rho T_g&=\rho \kappa_{(P4,P^3)}e^{3u}\;\;&\text{on}\;\;\partial M;\\
\frac{\partial u}{\partial n_g}&=0\;\;&\text{on}\;\;\partial M.
\end{split}
\right.
\end{equation}
\end{pro}
Now we are ready to make the proof of Theorem\;$\ref{eq:theorem3}$.\\
\vspace{10pt}

\noindent
\begin{pfn}{\sc of Theorem\;$\ref{eq:theorem3}$}\\
\vspace{2pt}

\noindent
By\;$\eqref{eq:cra}$\;and\;$\eqref{eq:pss}$\; there exists a sequence \;$\rho_{l}\rightarrow 1$\;and \;$u_{l}$\;such that the following holds :
\begin{equation}\nonumber
\left\{
\begin{split}
P^4_gu_l+2\rho_l Q_g&=0\;\;&\text{in}\;\;M;\\
P^3_gu_l+\rho_lT_g&=\rho \kappa_{(P4,P^3)}e^{3u_l};\;&\text{on}\;\;\partial M;\\
\frac{\partial u_l}{\partial n_g}&=0\;\;&\text{on}\;\;\partial M.
\end{split}
\right.
\end{equation}
Now  since \;$\kappa _{(P^4,P^3)}=\int_{M }Q_{g}dV_{g}+\int_{\partial M}T_gdS_g$\; then applying corollary\;$\ref{eq:compco}$\;with \;$Q_l=\rho_lQ_{g}$,\;$T_l=\rho_lT_g$\;and \;$\bar T_l=\rho_l \kappa_{(P^4,P^3)}$\;we have that\; $u_l$ is bounded in  \;$C^{4+\alpha}$\;for every \;$\alpha\in (0,1)$. Hence up to a subsequence it converges in \;$C^1(M)$\; to a solution of \;$\eqref{eq:bvps}$. Hence Theorem\;$\ref{eq:theorem3}$\;is proved.
\end{pfn}
\vspace{2pt}

\noindent
\begin{remark}
As said in the introduction, we now discuss how to settle the general case.\\
First of all, to deal with the remaining cases of {\em situation 1}, we proceed as in \cite{dm}. To obtain Moser-Trudinger type inequality and its improvement we impose the additional condition \;$\|\hat u\|\leq C$\;where\;$\hat u$\;is the component of \;$u$\; in the direct sum of the negative eigenspaces. Furthermore another aspect has to be considered, that is not only \;$e^{3u}$\;can concentrate but also \;$\|\hat u\|$\;can also tend to infinity. And to deal with this we have to substitute the set\;$\partial M_k$\;with an other one,\;$A_{k,\bar k}$\; which is defined in terms of the integer \;$k$ \;(given in\;$\eqref{eq:range}$)\;and the number \;$\bar k$\;of negative eigenvalues of \;$P^{4,3}_g$, as is done in \;$\cite{dm}$. This also requires suitable adaptation of the min-max scheme and of the monotonicity formula in Lemma\;$\ref{eq:noinc}$, which in general becomes
\begin{equation}\nonumber
\rho\rightarrow \frac{\ov{II}_{\rho}}{\rho}-C\rho\;\;\;\text{is non-increasing in }\;\;[1-\rho_{0},1+\rho_{0}];
\end{equation}
for a fixed constant \;$C>0$.\\
As already mentioned in the introduction, see Remark, to treat the {\em situation 1}, we only need to consider the case\;$\bar k\neq 0$. In this  case the same arguments as in \cite{dm} apply without any modifications.
\end{remark}

\footnotetext[1]{E-mail addresses:  ndiaye@sissa.it}


\begin{thebibliography}{99}



\bibitem {au} Aubin T.,{\em Nonlinear Analysis on manifold, Monge-Ampere equations}, Springer-Verlag, 1982.
\bibitem{bran1} Branson T.P., {\em The functional determinant}, Global Analysis Research Center Lecture Note Series, Number 4, Seoul National University (1993).

\bibitem{bran2} Branson T.P., {\em Differential operators canonically associated to a conformal structure}, Math. scand., 57-2 (1995), 293-345.

\bibitem{bo}Branson T.P., Oersted., {\em Explicit functional determinants in four dimensions}, Proc. Amer. Math. Soc 113-3(1991), 669-682.

\bibitem{br} Brendle S., {\em A family of curvature flows on surfaces with boundary}, Math. Z. 241, 829-869 (2002).

\bibitem{bren} Brendle S.,{\em Global existence and convergence for a higher order flow in conformal geometry}, Ann. of Math. 158 (2003),323-343.

\bibitem{bm} Brezis H., Merle F., {\em Uniform estimates and blow-up
behavior for solutions of $-\Delta u =V(x) e\sp u$ in two
dimensions} Commun. Partial Differ. Equations 16-8/9 (1991),
1223-1253.

\bibitem{cq1} Chang S.Y.A., Qing J.,{\em The Zeta Functional Determinants on manifolds with boundary 1. The Formula}, Journal of Functional Analysis 147, 327-362 (1997)

\bibitem{cq2} Chang S.Y.A., Qing J.,{\em The Zeta Functional Determinants on manifolds with boundary II. Extremal Metrics and Compactness of Isospectral Set}, Journal of Functional Analysis 147, 363-399 (1997)

\bibitem{cqy1} Chang S.Y.A., Qing J.,., Yang P.C.,{\em Compactification of a class of conformally flat 4-manifold}, Invent. Math. 142-1(2000), 65-93.

\bibitem{cy}  Chang S.Y.A.,  Yang P.C.,{\em Extremal metrics of zeta functional determinants on 4-manifolds}, ann. of Math. 142(1995), 171-212.
\bibitem{cy1}  Chang S.Y.A.,  Yang P.C.,{\em On a fourth order curvature invariant}.

\bibitem{cs} Chen, S.S., {\em Conformal deformation on manifolds with boundary}, Preprint
%


\bibitem{dj} Djadli Z., {\em Existence result for the mean field problem
on Riemann surfaces of all genus}, preprint.


\bibitem{djlw} Ding W., Jost J., Li J., Wang G., {\em Existence results
for mean field equations}, Ann. Inst. Henri Poincar\'e, Anal. Non
Lin�ire 16-5 (1999), 653-666.


\bibitem{dm1} Djadli Z., Malchiodi A., {\em A fourth order uniformization
theorem on some four manifolds with large total $Q$-curvature},
C.R.A.S., 340 (2005), 341-346.

\bibitem{dm} Djadli Z., Malchiodi A., {\em Existence of
conformal metrics with constant $Q$-curvature}, Ann. of Math, to appear.

\bibitem{dr} Druet O., Robert F., {\em Bubbling phenomena for fourth-order four-dimensional PDEs with exponential growth}, Proc. Amer. Math. Soc 134(2006) no.3, 897-908

\bibitem{es} Escobar J.F., {\em Conformal deformation of a Riemannian metric to a scalar flat metric with constant mean curvature on the boundary}, Ann. of Math. (2) 136 (1992), no. 1, 1--50 

\bibitem{fg} Fefferman C., Graham C,R., {\em Q-curvature and Poincar\'e metrics}, Mathematical Research Letters 9, 139-151(2002).




\bibitem{fg1} Fefferman C., Graham C., {\em Conformal invariants}, In Elie Cartan et les mathematiques d'aujourd'hui. Asterisque (1985), 95-116.

\bibitem{gt} Gilbar D., Trudinger N., {\em  Elliptic Partial Differential Equations of Second Order}, 2nd edition, Springr-Verlag, 1983.

\bibitem{gjms} Graham C,R., Jenne R., Mason L., Sparling G., {\em Conformally invariant powers of the laplacian, I:existence},  J.London Math.Soc 46(1992), no.2, 557-565.

\bibitem{gz} Graham C,R.,Zworsky M., {\em Scattering matrix in conformal geometry}. Invent math. 152,89-118(2003).


%

%








\bibitem{lm} Lions J.L., Magenes E., {\em Non-Homogeneous Boundary Value Problems and Applications}, Springer-Verlag Berlin Heidelberg New York 1972 (Volume 1).








\bibitem{mn} Malchiodi A., Ndiaye C.B {\em Some existence results for the Toda system on closed surface}, Rend. Mat. Acc. lincei, to appear.


\bibitem{nd} Ndiaye C.B., {\em Constant Q-curvature metrics in arbitrary dimension}, J. Funct. Anal
, to appear.


\bibitem{nd1} Ndiaye C.B., {\em Conformal metrics with constant \;$Q$-curvarure for manifolds with boundary}, preprint 2007.


\bibitem{nd2} Ndiaye C.B., {\em Curvature flows on four manifolds with boundary}, preprint 2007.


\bibitem{p1} Paneitz S., {\em A quartic conformally covariant differential operator for arbitrary pseudo-Riemannian manifolds}, preprint, 1983.


\bibitem{p2} Paneitz S., {\em Essential unitarization of symplectics and applications to field quantization}, J. Funct. Anal. 48-3 (1982), 310-359.





%
%








\bibitem{s} Stein E.M., Weiss G., {\em Introduction to Fourier Analysis on Euclidean spaces}, Princeton, New Jersey, Princeton University Press.

\bibitem{str} Struwe M., {\em The existence of surfaces of constant mean
curvature with free boundaries}, Acta Math. 160 -1/2(1988), 19-64.


%
%






\bibitem{xu} Xu Xingwang., {\em Uniqueness and non-existence theorems for conformally invariant equations},   Journal of Functional Analysis. 222(2005) 1-28. 

\end{thebibliography}
\end{document}